\documentclass[11pt]{article}
\usepackage{amsmath,amssymb,amsfonts}
\usepackage{algorithmic}
\usepackage{graphicx}
\usepackage{textcomp}

\usepackage{subcaption}

\newtheorem{lemma}{Lemma}
\newtheorem{theorem}{Theorem}
\newtheorem{example}{Example}

\title{Time and frequency domain low order, low frequency  approximation of mechanical systems}
\author{ Hans Zwart\footnotemark[1]\,  Dani\"el W.M.\  Veldman\footnotemark[2]\,  and Sahar F.\ Sharifi\footnotemark[3]}

\begin{document}
\maketitle

\renewcommand{\thefootnote}{\fnsymbol{footnote}}

\footnotetext[1]{H.\ Zwart is with the Department of Applied Mathematics, Faculty of Electrical Engineering, Mathematics and Computer Science, University of Twente, P.O. Box 217, 7500 AE, Enschede, The Netherlands and with the Department of Mechanical Engineering, Eindhoven University of Technology, P.O.\ Box 513, 5600 MB Eindhoven, The Netherlands (e-mail: H.J.Zwart@utwente.nl).}
\footnotetext[2]{D.W.M.\ Veldman is with the Department of Data Science, Friedrich-Alexander Universit\"at Erlangen-N\"urnberg, 91052 Erlangen, Bavaria, Germany (e-mail: d.w.m.veldman@gmail.com)}
\footnotetext[3]{S.F.\ Sharifi was with the Department of Applied Mathematics, Faculty of Electrical Engineering, Mathematics and Computer Science, University of Twente, P.O. Box 217, 7500 AE, Enschede, The Netherlands. She is now with TEMPRESS, Radeweg 31, 8171 MD Vaassen, The Netherlands (e-mail: ssharifi@tempress.nl). The research has received funding from the European Union’s Horizon 2020 research and
innovation programme under the Marie Sklodowska-Curie grant agreement No.\
765579.}

%\maketitle

\begin{abstract}
Control design for linear, time-invariant mechanical systems typically requires an accurate low-order approximation in the low frequency range. For example a series expansion of the transfer function around zero consisting of a mass, velocity, and compliance term. Because computing such a series expansion of the transfer function can be cumbersome, a new method to compute low-order approximations of mechanical systems is developed in this paper. The method does not require an explicit expression for the transfer function, which is not always available for infinite-dimensional systems. The advantages of the proposed method is demonstrated in three examples. 
\end{abstract}

\noindent
{\bf Keywords}:
Compliance, Model Order Reduction, Partial Differential Equations, Boundary control systems. Transfer functions

\noindent
{\bf AMS subject classification.} 93B11, 93C05, 93C20, 93C80,  93C95.

\section{Introduction}
\label{sec:1}
For control design of linear, time-invariant (LTI) systems,  the low frequency range is often the most important and an accurate model is needed in this range. When the system has no poles on the imaginary axis, then the steady-state gain, i.e, the transfer function at frequency zero, serves as an initial approximation. However, for mechanical systems, with a force input and a position output, there will always be a double pole in the origin. This is a direct consequence of Newton law, i.e., force is proportional to acceleration. For these systems an initial approximation of the transfer function would be $1/(ms^2)$, with $m$ the (total) mass. However, this approximation might capture too little of the system. This is the case in e.g.\ the design of feedforward tracking control, see Kontaras et al \cite{Nikos}. Then also terms proportional to $1/s$ and the constant should be taken into account. The constant term is known as the compliance. 

For finite-dimensional systems, the compliance is typically computed based on a complete eigenvalue decomposition, which can be computationally demanding when the system has a large number of degrees of freedom. For infinite-dimensional systems, finding a low-order approximation typically requires determining the transfer function first. However, finding an explicit expression for the transfer function is typically impossible when the model is given as a partial differential equation on domains with more than one spatial dimension or with spatially dependent coefficients. A method to determine the compliance in an easier way, i.e.\ without computing a full basis of eigenvectors or finding a closed-form analytic expression for the transfer function, is thus needed. 

In this paper, a time-domain characterization of the first three terms that appear in the series expansion of the transfer function of a mechanical system around zero (i.e.\ the mass, velocity, and compliance terms) is provided. This enables us to determine the compliance from a set of conditions in the time-domain that are applicable to partial differential equations with spatially dependent coefficients on general spatial domains. The derived characterization also shows that the compliance can be determined from a time-domain simulation of the mechanical system. Since many efficient methods for time-domain simulations exist, this provides a quick way of obtaining a low order approximation of the mechanical system at hand. 

The remainder of this paper is organized as follows. In the following section, the frequency domain definition of the compliance is given together with a sketch of the main idea of the proposed time-domain approach and how it can be used to determine the compliance from a time-domain simulation. The equivalence between the time-domain and frequency-domain characterization of the compliance are proven in Sections \ref{sec:3} and \ref{sec:4}. Section \ref{sec:3} does this for finite-dimensional systems and infinite-dimensional systems with internal control and Section \ref{sec:4} for boundary control systems. In Section \ref{sec:5}, the advantages of the proposed time-domain approach to the compliance are illustrated in three examples: a vibrating string with structural damping, a Euler-Bernoulli beam with a viscous damper at one of its endpoints, and a vibrating plate with structural damping. In the latter example, the compliance is also determined from a time-domain simulation. 

\section{Compliance, frequency and time domain}
\label{sec:2}

When a mechanical system contains rigid body modes, the transfer function $G(s)$ from a force input to a position output typically includes a double-integrator-part. This mean that $s^2G(s)$ is regular near $s = 0$ and may be expanded in a Taylor series. This leads to the approximation
\begin{align}
  G(s) &= \frac{1}{s^2} \left( G_2 + G_1 s + G_0s^2 \right) + G_{st}(s) 
  \nonumber \\ 
  &= \frac{G_2}{s^2} + \frac{G_1}{s} + G_0 + G_{st}(s), 
\label{eq:Gs}
\end{align}
where $G_{2}^{-1}$, $G_{1}$, and $G_0$ denote the mass, velocity, and compliance terms in the expansion of the transfer function, respectively. $G_{st}(s)$ is the stable part, i.e, no poles in the closed right half plane, and captures the remainder of order $O(s)$. The mass and velocity part $G_2$ and $G_1$ are typically related to the rigid body modes of the systems and $G_{flex}(s) = G_0 + G_{st}(s)$ captures the flexible part related to the (possibly infinitely many) non-rigid modes. 

When the transfer function $G(s)$ is available, it is clear from \eqref{eq:Gs} that $G_2$, $G_1$, and $G_0$ may be computed as
\begin{align}
G_2 &= \lim_{s\rightarrow 0} s^2G(s),  \label{eq:G2}\\
G_1 &= \lim_{s\rightarrow 0} s \left(G(s) - \frac{G_2}{s^2} \right),  \label{eq:G1} \\
G_0 &= \lim_{s\rightarrow 0} \left(G(s) - \frac{G_2}{s^2}-\frac{G_1}{s}\right).  \label{eq:G0}
\end{align}
For infinite-dimensional systems with constant coefficients on a one-dimensional domain, it is often possible to determine a closed-form analytic expression for the transfer function, see, e.g., \cite{Ruth} or\cite{curtain2009}, and a low frequency approximation can be determined from \eqref{eq:G2}--\eqref{eq:G0}. However, in many other situations, a closed-form analytic expression for the transfer function is not available or computationally demanding to evaluate. Determining a low-order approximation from \eqref{eq:G2}--\eqref{eq:G0} is then impossible. Even when a closed-form analytic expression is available, determining $G_2$, $G_1$, and $G_0$ from \eqref{eq:G2}--\eqref{eq:G0} can be very cumbersome, as will be illustrated in Subsection \ref{Ex-beam}. 

Motivated by these difficulties, a time-domain approach based on the step response of the system, described by the transfer function $G(s)$, to determine $G_2$, $G_1$, and $G_0$, is proposed. In particular, let $y(t)$ denote the output of the system described by $G(s)$ resulting from the step input $u(t)=u_0$ for $t > 0$ and zero initial conditions. Taking the inverse Laplace transform in \eqref{eq:Gs} readily shows that
\begin{equation}
\label{eq-sol}
    y(t)=G_2u_0\frac{t^2}{2} + G_1u_0 t + G_0u_0 +\omega_{st}(t),
\end{equation}
where $\omega_{st}(t)$ is the inverse Laplace transform of $G_{st}(s)\frac{u_0}{s}$. Since $G_{st}$ is stable and $G_{st}(0)=0$, $\omega_{st}(t)$ converges to zero as $t\rightarrow \infty$. The compliance function can thus be determined from the static part of the step response. 

This observation can now be used to determine $G_2$, $G_1$, and $G_0$ based on the step response. Although the above seems to indicate that this is straightforward, some care should be taken. This is due to the fact that the inverse Laplace transform assumes that the initial conditions are taken on $0^-$, i.e., on $\lim_{t \uparrow 0}$, whereas calculation in time domain take the initial condition at $0^+$. When the solution is continuous at $t=0$, then there is no difference, but this need not to hold for partial differential equations with boundary control, or when the control also enters with its derivative. 

Based on (\ref{eq-sol}), it is assumed that the response resulting from a step input $u(t)=u_0$ for $t > 0$ is of the form
\begin{equation}
\label{eq:omega}
  \omega(t) =\omega_2\frac{t^2}{2} + \omega_1t + \omega_0 +\omega_{st}(t),
\end{equation}
where $\omega_{st}(t)$ is stable, i.e.\ $\omega_{st}(t) \rightarrow 0$ for $t \rightarrow \infty$. 
Comparing \eqref{eq-sol} to \eqref{eq:omega}, it follows that 
\begin{equation}
G_2u_0 = \omega_2, \qquad
G_1u_0 = \omega_1, \qquad
G_0u_0 = \omega_0. \label{eq:Gk}
\end{equation}
Therefore, $G_2$, $G_1$, and $G_0$ could be determined from a set of step responses when the $u_0$'s span the input space. 

Finally, observe that the $\omega_2$, $\omega_1$, and $\omega_0$ in \eqref{eq:omega} can be determined from snapshots of $\omega(t)$. In particular, let $t_1 < t_2 < \ldots < t_n$ be time instances that are so large that $\omega_{st}(t_i)$ with $i \in \{1,2, \ldots, n \}$ is small/negligiable. Then \eqref{eq:omega} shows that
\begin{equation}
\begin{bmatrix}
\omega(t_1) \\ \omega(t_2) \\ \vdots \\ \omega(t_n)
\end{bmatrix} \approx \begin{bmatrix}
\tfrac{1}{2}t_1^2 & t_1 & 1 \\
\tfrac{1}{2}t_2^2 & t_2 & 1 \\
\vdots & \vdots & \vdots \\
\tfrac{1}{2}t_n^2 & t_n & 1
\end{bmatrix} \begin{bmatrix}
\omega_2 \\ \omega_1 \\ \omega_0
\end{bmatrix} . \label{eq:omegai}
\end{equation}
Estimates for $\omega_2$, $\omega_1$, and $\omega_0$ can thus be determined by minimizing the residue in \eqref{eq:omegai}. The resulting estimates $\tilde{\omega}_2$, $\tilde{\omega}_1$, and $\tilde{\omega}_0$ are given by
\begin{equation}
\begin{bmatrix}
\tilde{\omega}_2 \\ \tilde{\omega}_1 \\ \tilde{\omega}_0
\end{bmatrix} = \left( \sum_{i=1}^n \begin{bmatrix}
\tfrac{1}{2}t_i^2 \\ t_i \\ 1
\end{bmatrix} \begin{bmatrix}
\tfrac{1}{2}t_i^2 & t_i & 1
\end{bmatrix} \right)^{-1} \sum_{i=1}^n \begin{bmatrix}
\tfrac{1}{2}t_i^2 \\ t_i \\ 1
\end{bmatrix} \omega(t_i). \label{eq:omegai_sol}
\end{equation}
Note that there is often no need to implement \eqref{eq:omegai_sol} directly, but that the least squares solution to \eqref{eq:omegai} can be computed by built in commands, e.g.\ the \texttt{/}-operator in Matlab$^\copyright$.  The remaining residue gives an indication whether $\omega_{st}(t)$ is indeed negligible at the considered time instances $t_i$ or not.

%%%%%%%%%%%%%%%%%%%%%

\section{Compliances for internally controlled systems}
\label{sec:3}

In this section, the relation between the frequency domain and time domain characterization of the compliance function is studied for systems with internal control. Consider the following class of models:
\begin{align}
\begin{split}
  &\ddot{\omega}(t)+\mathcal{D}_{0}\dot{\omega}(t)+\mathcal{A}_{0}\omega(t) = B_0u(t) + B_1\dot{u}(t), \\
  &\omega(0)=0, \quad  \dot{\omega}(0)=0, 
\label{eqm-1}
 \end{split}
 \end{align}
 where $\mathcal{A}_{0}: \mathbf{D}(\mathcal{A}_{0}) \subseteq X\rightarrow X$, $\mathcal{D}_{0}: \mathbf{D}(\mathcal{D}_{0}) \subseteq X\rightarrow X$  are linear (possibly unbounded) operators on the Hilbert space $X$, $B_0$ and $B_1$ are bounded operators from another Hilbert space $U$ to the state space $X$. The operator ${\mathcal D}_0$ will represent the damping. Note that all physical parameters are hidden in the operators $\mathcal{A}_{0}$, ${\mathcal D}_0$, $B_0$, and $B_1$. The variable $u(t)$ is the input function/force, and $\omega(t)$ is the output/position. Thus $\dot{\omega}$ and $\ddot{\omega}$ denote the velocity and acceleration, respectively, and it is assumed that the damping depends on the velocity. The term $B_1\dot{u}(t)$ naturally appears in finite element models of certain partial differential equations with boundary control. This will be further illustrated in Section \ref{sec:4} and the example in Subsection \ref{example}.
 
Recall that $\omega$ is a classical solution of (\ref{eqm-1}), when it is twice differentiable, takes values in $\mathbf{D}(\mathcal{A}_{0})$, $\dot{\omega}$ takes its values in $\mathbf{D}(\mathcal{D}_{0})$, and $\omega$ satisfies (\ref{eqm-1}). 
 Our investigations are focused on the step response obtained by setting $u(t)=u_0$  for $t > 0$. When $B_1 \neq 0$, the step response will typically contain a jump in the velocity and, after the jump, the weak solution of (\ref{eqm-1}) equals the weak solution of
\begin{align}
\begin{split}
\label{eqm-2}
&\ddot{\omega}(t)+\mathcal{D}_{0}\dot{\omega}(t)+\mathcal{A}_{0}\omega(t)=B_0u_0,\\
&\omega(0)=0, \quad  \dot{\omega}(0^+)=B_1 u_0.
\end{split}
\end{align}
It is now assumed that this equation possesses a classical solution.  In particular, it is required that $B_1 u_0 \in \mathbf{D}(\mathcal{D}_{0})$. 
 
When $\mathcal{A}_0$ is selfadjoint, $\mathbf{D}(\mathcal{D}_0) = X$, and $-\mathcal{D}_0$ is dissipative, the classical solution of \eqref{eqm-2} on $[0,\infty)$ exists when $B_1u_0 \in \mathbf{D}(\mathcal{A}_{0}^{\frac{1}{2}})$, see, e.g., \cite[Example 2.3.5]{Ruth}.
%
%\textcolor{blue}{It is easy to see that a necessary condition for $\omega(t)$ being of the form \eqref{eq:omega} is that the step response does not contain any undamped oscillations, i.e.\ that all eigenfunctions of $\mathcal{A}_0$ that are not in the kernel of $\mathcal{A}_0$, also do not lie in the kernel of $\mathcal{D}_0$. When the eigenfunctions of $\mathcal{A}_0$ and $\mathcal{D}_0$ are the same (i.e.\ for proportional damping, see e.g.\ \cite{rixen2014}), this is also a sufficient condition which reduces to $\ker(\mathcal{D}_0) \subseteq \ker(\mathcal{A}_0)$. %The latter case covers a large class of engineering problems. 
%More generally, it can be shown that solutions of \eqref{eqm-2} decay when the flow generated by \eqref{eqm-2} with $\mathcal{D}_0 = 0$ is observable through $\mathcal{D}_0$, see, e.g., \cite{haraux1989, chill2019}. } 
%
\begin{theorem}
\label{Th-1}
The (classical) solution of \eqref{eqm-2} is of the form (\ref{eq:omega}) if and only if the following equations are satisfied
\begin{align}
  \label{eq:8a} 0=&\ {\mathcal A}_0\omega_2\\
   \label{eq:8b} 0=&\ {\mathcal D}_0\omega_2 + {\mathcal A}_0\omega_1\\
   \label{eq:8c} B_0u_0 =&\ \omega_2 + {\mathcal D}_0\omega_1 + {\mathcal A}_0\omega_0\\
    0 =&\ \ddot{\omega}_{st}(t) + {\mathcal D}_0 \dot{\omega}_{st}(t) + {\mathcal A}_0 \omega_{st}(t), \nonumber \\
    &\omega_{st}(0)=-\omega_0, \quad \dot{\omega}_{st}(0^+)=B_1u_0-\omega_1, \label{eq:8d}
\end{align}
and the (classical) solution $\omega_{st}(t)$ of \eqref{eq:8d} is stable, i.e.\ $\lim_{t\rightarrow\infty} \omega_{st}(t) = 0$. 
Furthermore, if for every $u_0 \in U$ such solution exists, then the transfer function is of the form \eqref{eq:Gs} and \eqref{eq:Gk} holds. 
\end{theorem}
{\bf Proof} For the first part just substitute the assumed solution \eqref{eq:omega} into the differential equation \eqref{eqm-2}, and equate the terms $t^2$, $t$, etc. Since $\omega(0)$ must be zero, it follows from \eqref{eq:omega} that $\omega_{st}(0)$ must be $-\omega_0$. Similarly, since $\dot{\omega}(0^+)= B_1 u_0$, it follows that $\dot{\omega}_{st}(0^+) =B_1u_0 -\omega_1$. 

For the second part, take the Laplace transform of the proposed solution (\ref{eq:omega}). Standard Laplace theory gives that this equals
\[
  Y(s) = \frac{\omega_2}{s^3} + \frac{\omega_1}{s^2} + \frac{\omega_0}{s} + \Omega_{st}(s).
\]
Also, taking the Laplace transform in \eqref{eqm-1}, the following relation between the Laplace transform of any solution $\Omega(s)$ starting from zero initial conditions and the Laplace transform of the input $U(s)$ is obtained
 \[
   (s^2I + s{\mathcal D}_0+ {\mathcal A}_0)\Omega(s) = (B_0 + sB_1)U(s).
 \]
If it can be shown that if $Y(s)$ satisfies the above equation for $\Omega(s)$ with $U(s) = u_0 / s$, then, by uniqueness and by definition of transfer function $G(s)$, see \cite[Section 10.2]{Ruth}, it follows that $Y(s)=\Omega(s) = G(s)u_0/s$. Because of the expression for $Y(s)$ given above, it then follows that $G(s)$ is of the form \eqref{eq:Gs} and that \eqref{eq:Gk} holds. 

It thus remains to show that $Y(s)$ satisfies the equation for $\Omega(s)$ above with $U(s) = u_0/s$. To this end, note that
 \begin{align*}
   (s^2I + s&{\mathcal D}_0 + {\mathcal A}_0)Y(s) \\
   =&\  (s^2I + s{\mathcal D}_0 + {\mathcal A}_0) \left(\frac{\omega_2}{s^3} + \frac{\omega_1}{s^2} + \frac{\omega_0}{s} + \Omega_{st}(s)\right)\\
     =&\  \frac{{\mathcal A}_0\omega_2}{s^3} + \left[ {\mathcal D}_0\omega_2 + {\mathcal A}_0\omega_1\right] \frac{1}{s^2} +  \left[ \omega_2 + D_0\omega_1 +A_0\omega_0\right] \frac{1}{s} \\
     &\ + (s^2I + s{\mathcal D}_0 + {\mathcal A}_0)\Omega_{st}(s) + \left[ \omega_1 + {\mathcal D}_0\omega_0\right] + s \omega_0\\
     =&\ 0+ \frac{B_0u_0}{s} 
      + (s^2I + s{\mathcal D}_0 + {\mathcal A}_0)\Omega_{st}(s)+ \left[ \omega_1 + {\mathcal D}_0\omega_0\right] + s \omega_0,
\end{align*}
where (\ref{eq:8a})--(\ref{eq:8c}) have been used. 

Because $\omega_{st}(t)$ is a solution of the differential equation (\ref{eq:8d}), its Laplace transform will satisfy a corresponding algebraic equation. However, some care should be taken. Normally the Laplace transform of $f(t)$ is defined as $F(s)=\int_{0^-}^\infty f(t) e^{-st} dt$. However, when $f$ is continuous at zero, then this integral equals $\int_{0^+}^\infty f(t) e^{-st} dt$. So $F(s)$ does not change, but the differentiation rule changes. Namely, the Laplace transform of $\dot{f}(t)$ equals $sF(s) - f(0^+)$. Similarly for higher order derivatives. Since (\ref{eq:8d}) has a condition on $0^+$, the integral $\int_{0^+}^\infty f(t) e^{-st} dt$ is used for the Laplace transform. Note that since $\omega_{st}$ is continuous at zero, its Laplace transform does not change. Knowing this, the Laplace transform of the differential equation for $\omega_{st}$ in \eqref{eq:8d} yields
\begin{align*}
  0=&\ s^2 \Omega_{st}(s) - s \omega_{st}(0) - \dot\omega_{st}(0^+) 
  + {\mathcal D}_0 \left[ s \Omega_{st}(s) - \omega_{st}(0) \right] + {\mathcal A}_0\Omega_{st}(s) \\
  =&\ (s^2I + s{\mathcal D}_0 + {\mathcal A}_0)\Omega_{st}(s) + s\omega_0 - B_1u_0 + \omega_1 + {\mathcal D}_0 \omega_0. 
\end{align*}
Subtracting this from the expression for $(s^2I + s{\mathcal D}_0 + {\mathcal A}_0)Y(s)$ above yields
\[ 
(s^2I + s{\mathcal D}_0 + {\mathcal A}_0)Y(s) = \frac{B_0u_0}{s} + B_1u_0 = \left(B_0 + sB_1 \right) \frac{u_0}{s}. 
\]
It thus follows that $\Omega(s) = Y(s)$ when $U(s) = u_0/s$, which completes the proof. 
 \hfill$\Box$
 \medskip

The above theorem gives conditions for finding the $\omega_k$'s. With additional structure on the operators more can be said about these solutions. Recall that the operator $Q$ with domain $\mathbf{D}(Q)$ is {\em symmetric}\/ when $\langle Qz_1,z_2\rangle = \langle z_1,Qz_2\rangle$ for all $z_1, z_2 \in {\mathbf D}(Q)$. Furthermore, it is non-negative when $\langle Qz,z\rangle \geq 0$ for $z\in {\mathbf D}(Q)$. 
\begin{theorem} \label{thm:om1}
When ${\mathcal A}_0, {\mathcal D}_0$ are symmetric and ${\mathcal D}_0$ is non-negative, then the set of equalities (\ref{eq:8a})--(\ref{eq:8c}) imply that ${\mathcal D}_0\omega_2=0$ and ${\mathcal A}_0\omega_1=0$. Moreover, equation \eqref{eq:8d} combined with the limit behaviour of its solution implies that for every $v \in \ker(\mathcal{A}_0) \cap \ker(\mathcal{D}_0)$
\begin{equation}
\langle v, \omega_0 \rangle = 0, \qquad  \langle v, \omega_1 \rangle = \langle v, B_1u_0 \rangle.  \label{eq:thm2}
\end{equation} 
Furthermore, if $\ker(\mathcal A_0) \subseteq \ker(\mathcal{D}_0)$ and $B_1 u_0 \in \ker(\mathcal{A}_0)^\perp$, then $\omega_1 = 0$. 
\end{theorem}
{\bf Proof}
Taking the inner product of \eqref{eq:8b} with $\omega_2$ shows
\[
  0= \langle \omega_2, {\mathcal D}_0\omega_2 + {\mathcal A}_0\omega_1\rangle =  \langle \omega_2, {\mathcal D}_0\omega_2 \rangle + \langle {\mathcal A}_0\omega_2,\omega_1\rangle 
  = \langle \omega_2, {\mathcal D}_0\omega_2 \rangle,
\]
where it was used that ${\mathcal A}_0$ is symmetric and that ${\mathcal A}_0\omega_2 =0$, see (\ref{eq:8a}). Since ${\mathcal D}_0$ is non-negative this equation implies that ${\mathcal D}_0\omega_2=0$.
% \textcolor{blue}{Proof: Let $q \in D({\mathcal D_0})$ and $t \in {\mathbb R}$, then
%\[
%  0\leq \langle \omega_2 + t q, {\mathcal D}_0(\omega_2 +t q) \rangle = 0 + 2 t \langle q, {\mathcal D}_0\omega_2 \rangle + t^2  \langle  q, {\mathcal D}_0(q) \rangle.
%\]
%This can only hold when $\langle q, {\mathcal D}_0\omega_2 \rangle =0$. Since $q$ is arbitrary in the dense set $D({\mathcal D_0})$ it follows that  ${\mathcal D}_0\omega_2 =0$
%}. 
That $\mathcal{A}_0 \omega_1 = 0$ now follows from \eqref{eq:8b}. 

For \eqref{eq:thm2}, let $v \in \ker(\mathcal{A}_0) \cap \ker(\mathcal{D}_0)$ and set
\[
  v_{st}(t) = \langle v, \omega_{st}(t) \rangle.
\]
By the differential equation of $\omega_{st}$ in \eqref{eq:8d}
\begin{align*}
 \ddot{v}_{st}(t) =&\ -  \langle v, {\mathcal D}_0 \dot{\omega}_{st}(t) \rangle -  \langle v, {\mathcal A}_0 \omega_{st}(t) \rangle\\
=&\ -\langle {\mathcal D}_0v, \dot{\omega}_{st}(t) \rangle -  \langle {\mathcal A}_0 v, \omega_{st}(t) \rangle = 0,
\end{align*} 
where it was used that ${\mathcal A}_0$, ${\mathcal D}_0$ are symmetric and that $v \in \ker(\mathcal{A}_0) \cap \ker(\mathcal{D}_0)$. 
Therefore, 
\[
  v_{st}(t)  = \alpha + \beta t.
\]
Since $\omega_{st}(t)$ is stable, it must hold that $\alpha=\beta=0$, and thus $v_{st}(t)$ is identically zero. By definition of $v_{st}(t)$, this implies that $\langle v, \omega_{st}(0) \rangle = 0$ and $\langle v, \dot{\omega}_{st}(0^+) \rangle = 0$. The initial conditions for $\omega_{st}(t)$ in \eqref{eq:8d} now give \eqref{eq:thm2}. 

The last assertion follows by taking $v = \omega_1$ in the second condition in \eqref{eq:thm2}, which is possible because it has been shown before that $\omega_1 \in \ker(\mathcal{A}_0)$. 
\hfill$\Box$

The following example shows that $\omega_1$ can be unequal to zero when $B_1 = 0$.
\begin{example}
\label{E:om1}
Consider the following system on ${\mathbb R}^3$
\[
  \ddot{x}(t) + {\mathcal D}_0 \dot{x}(t) + {\mathcal A}_0 x(t) = B_0u_0, \quad x(0)=0, \dot{x}(0)= 0,
\]
with
\[
  {\mathcal D}_0 = \begin{bmatrix} 0 & 0 & 0 \\ 0 & 1 & 0 \\ 0 & 0 & 1 
\end{bmatrix}, \quad  
{\mathcal A}_0 = \begin{bmatrix} 0 & 0 & 0 \\ 0 & 0 & 0 \\ 0 & 0 & 1 \end{bmatrix}, \quad
B_0 = \begin{bmatrix} 1 \\ 1 \\ 1 \end{bmatrix}.
 \]
Note that ${\mathcal A}_0$ and ${\mathcal D}_0$ are symmetric.
For $u_0=1$, the above vector valued ODE becomes
\begin{align*}
  \ddot{x}_1(t) = 1,\quad& x_1(0)=0, \dot{x}_1(0)= 0,\\
  \ddot{x}_2(t) + \dot{x}_2(t) = 1,\quad&x_2(0)=0, \dot{x}_2(0)= 0,\\
  \ddot{x}_3(t) + \dot{x}_3(t) + x_3(t) = 1,\quad & x_3(0)=0, \dot{x}_3(0)= 0.
 \end{align*}
 The solutions are $x_1(t) = \frac{1}{2} t^2$, $x_2(t)= t -1 + e^{-t} $, and 
 \[
  x_3(t) = -e^{-\frac{1}{2}t} \cos(\frac{\sqrt{3}}{2} t) - \frac{1}{3} \sqrt{3} e^{-\frac{1}{2}t} \sin(\frac{\sqrt{3}}{2} t) +1.
\]
Hence $\omega_1=\left[0,1,0\right]^T$, which is unequal to the zero vector.
The reason for $\omega_1 \neq 0$ is that $\ker({\mathcal A}_0) \not\subseteq \ker({\mathcal D}_0)$.
\end{example}

In the last part of this section, a formula for the compliance when ${\mathcal A}_0$ is a non-negative self-adjoint operator with compact resolvent is derived. Let $\lambda_n$ and $\varphi_n$ ($n \in \mathbb{N}$) denote the eigenvalues and eigenvectors of $\mathcal{A}_0$. Since $\mathcal{A}_{0}$ is self-adjoint and has compact resolvent, the eigenvectors $\varphi_n$ can be chosen as an orthonormal basis for $X$. 
%\begin{equation}
%\label{eq:2.4}
%   X=\overline{\mathrm{span}}\{\varphi_{1},\varphi_{2},\varphi_{3},\cdots,\varphi_{n+1}, \cdots\}. 
%\end{equation}
Furthermore, it is assumed that $\mathcal{D}_0$ is diagonalizable on the same basis of eigenvectors as $\mathcal{A}_0$, i.e.\ \emph{proportional damping}, see e.g.\ \cite{rixen2014}. This means that $\mathcal{A}_0$ and $\mathcal{D}_0$ can be represented as
\begin{equation}
\mathcal{A}_0 \omega = \sum_{n \in \mathbb{N}} \lambda_n \langle \varphi_n, \omega \rangle \varphi_n, \quad \mathcal{D}_0 \omega = \sum_{n \in \mathbb{N}} \mu_n \langle \varphi_n, \omega \rangle \varphi_n, \label{eq:A0D0_spectral}
\end{equation}
with $\lambda_n, \mu_n \geq 0$. 
\begin{theorem} 
\label{thrm1}
%\normalfont  
Let ${\mathcal A}_0$ and $\mathcal{D}_0$ be non-negative self-adjoint operators with compact resolvent satisfying $\ker(\mathcal{D}_0) \subseteq \ker(\mathcal{A}_0)$. Let the eigenvalues and normalized eigenfunctions of $\mathcal{A}_0$ be $\lambda_n$ and $\varphi_{n}$, $n \in {\mathbb N}$, respectively. Let $\mathcal{D}_0$ have the same eigenfunctions $\varphi_n$ and let its eigenvalues be $\mu_n$, $n \in \mathbb{N}$. 
Then the $\omega_0$, $\omega_1$, and $\omega_2$ satisfying \eqref{eq:8a}--\eqref{eq:8d} are
\begin{equation} 
\label{w22a}
   \omega_{2}=\sum_{n\in \mathbb{J}}\alpha_n \varphi_n, \quad 
   \omega_1= \sum_{n\in {\mathbb I}}\beta_n \varphi_n, \quad
   \omega_{0}=\sum_{n\in \mathbb{N} \backslash \mathbb{J}} \gamma_n \varphi_n,
\end{equation}
where $\mathbb{I} = \{ n\in {\mathbb N} \mid \lambda_n = 0 \}$, $\mathbb{J} = \{ n \in {\mathbb N} \mid \mu_n = 0 \}$, and
\begin{align}
\label{alfa+beta1}
        \alpha_{n}=&\ \langle\varphi_{n},B_0u_{0}\rangle, \qquad  
        \beta_{n}= \begin{cases}
         \langle\varphi_{n},B_1u_{0}\rangle & n \in \mathbb{J} \\
        \frac{\langle\varphi_{n},B_0u_{0}\rangle}{\mu_n}, & n \in \mathbb{I} \backslash \mathbb{J}
        \end{cases}, \\
        \label{gamma}
        \gamma_{n}=&\ \begin{cases}
        \frac{\langle \varphi_{n}, B_0u_{0} \rangle}{ \lambda_n} & n \in \mathbb{N} \backslash \mathbb{I} \\
        \frac{\langle \varphi_n, B_1 u_0 \rangle}{ \mu_n} - \frac{\langle \varphi_n, B_0u_0 \rangle}{\mu_n^2}  & n \in \mathbb{I} \backslash \mathbb{J}
        \end{cases} .
\end{align}
\end{theorem}
{\bf Proof} 
Theorem \ref{thm:om1} shows that $\omega_2 \in \ker(\mathcal{D}_0)$ and $\omega_1 \in \ker(\mathcal{A}_0)$ which implies that $\omega_2$ and $\omega_1$ are of the form in \eqref{w22a} for certain coefficients $\alpha_n$ and $\beta_n$. 

Note that $\ker(\mathcal{D}_0) \subseteq \ker(\mathcal{A}_0)$ implies that $\mathbb{J} \subseteq \mathbb{I}$. Projecting \eqref{eq:8c} on the basis functions $\varphi_n$ with $n \in \mathbb{J}$ thus yields the expression for $\alpha_n$ in \eqref{alfa+beta1}. Applying the same procedure for basis functions $\varphi_n$ with $n \in \mathbb{I} \backslash \mathbb{J}$ yields the expression for $\beta_n$ in \eqref{alfa+beta1} for $n \in \mathbb{I} \backslash \mathbb{J}$, and projecting \eqref{eq:8c} on $\varphi_n$ with $n \in \mathbb{N} \backslash \mathbb{I}$ yields the formula for $\gamma_n$ in \eqref{gamma} for $n \in \mathbb{N} \backslash \mathbb{I}$. 

The solution of the second order differential equation (\ref{eq:8d}) is of the form
\begin{align}
\nonumber
  \omega_{st}(t) =&\ \sum_{n \in \mathbb{J}} (\gamma_{n,1} + t \gamma_{n,2}) \varphi_n
  + 
  \sum_{n \in \mathbb{I} \backslash \mathbb{J}}  
(\gamma_{n,1} + e^{-\mu_n t} \gamma_{n,2})  \varphi_n
  \\
  \label{omega_st}
  &\ + \sum_{n\in \mathbb{N} \backslash \mathbb{J}} \left(\gamma_{n,1} e^{\mu_{n,1} t} + \gamma_{n,2} e^{\mu_{n,2} t} \right) \varphi_n, 
\end{align}
for certain coefficients $\gamma_{n,1}$ and $\gamma_{n,2}$ and with
\[ 
  \mu_{n,12} = 
  \frac{-\mu_n \pm \sqrt{\mu_n^2 - 4 \lambda_n}}{2}. 
  \]
In order to have $\lim_{t \rightarrow \infty} \omega_{st}(t) = 0$, $\gamma_{n,1}$ should be zero for $n \in \mathbb{I}$ and $\gamma_{n,2}$ should be zero for $n \in \mathbb{J}$. 
Therefore, the expression for $\omega_{st}(t)$ above shows that 
\begin{align*}
\langle \varphi_n, \omega_{st}(0) \rangle = \langle \varphi_n, \dot{\omega}_{st}(0^+) \rangle &= 0, \quad \qquad n \in \mathbb{J}, \\
\mu_n \langle \varphi_n, \omega_{st}(0) \rangle + \langle \varphi_n, \dot{\omega}_{st}(0^+) \rangle &= 0, \quad\qquad  n \in \mathbb{I} \backslash \mathbb{J}.
\end{align*}
Inserting the initial conditions in \eqref{eq:8d} into the first condition yields that $\omega_0$ is of the form in \eqref{w22a} and that $\beta_n$ is as in \eqref{alfa+beta1} for $n \in \mathbb{J}$. The condition for $n \in \mathbb{I} \backslash \mathbb{J}$ yields the expression for $\gamma_n$ in \eqref{gamma} with $n \in \mathbb{I} \backslash \mathbb{J}$. %\marginpar{I Get the minus of the expression. Indeed the sign was wrong. Fixed now. }
\hfill$\Box$
%\medskip

%%%%%%%%%%%%%%%%%%%%%%%%%%%%

\section{Compliances for boundary control systems}
\label{sec:4}

Consider the following system with control at the boundary
\begin{align} 
\begin{split}
\label{eq:bc_sys}
&\ddot{\omega}(t) + {\mathcal D} \dot{\omega}(t) + {\mathcal A} \omega(t) = 0, \\
&{\mathcal B}_0 \omega(t) + \mathcal{B}_1 \dot{\omega}(t) = u(t), \\
&\omega(0)=0,  \dot{\omega}(0)=0,
\end{split}
\end{align}
where $\mathcal{A}: \mathbf{D}(\mathcal{A}) \subseteq X\rightarrow X$, $\mathcal{D}: \mathbf{D}(\mathcal{A}) \subseteq \mathbf{D}(\mathcal{D}) \subseteq X\rightarrow X$ are linear operators on a Hilbert space $X$, and $\mathcal{B}_0 : \mathbf{D}(\mathcal{A})\subseteq \mathbf{D}(\mathcal{B}_0) \subseteq X \rightarrow U$ and $\mathcal{B}_1 : \mathbf{D}(\mathcal{A})\subseteq \mathbf{D}(\mathcal{B}_1) \subseteq X\rightarrow U$, where $U$ is another Hilbert space, and $u(t)$ represents the control action at the boundary. The existence theory of boundary control systems can for example be found in \cite[Chapter 10]{Ruth}. Furthermore, denote
\begin{align}
\label{eq:23}
  {\mathcal A}_0 =&\ {\mathcal A}, \qquad \mathbf{D}({\mathcal A}_0) = \mathbf{D}({\mathcal A}) \cap \ker{\mathcal B}_0,\\
\label{eq:24}
  {\mathcal D}_0 =&\ {\mathcal D}, \qquad \mathbf{D}({\mathcal D}_0) = \mathbf{D}({\mathcal D}) \cap \ker{\mathcal B}_0,
\end{align}
and assume there exist operators 
$\mathcal{C}_A : \mathbf{D}(\mathcal{A}_0^*) \rightarrow U$ and $\mathcal{C}_D : \mathbf{D}(\mathcal{D}_0^*) \rightarrow U$ such that for all $\phi \in \mathbf{D}(\mathcal{A}_0^*)$ and $\omega \in \mathbf{D}(\mathcal{A})$
\begin{align}
\label{eq:CA}
\langle \mathcal{A}_0^* \phi, \omega \rangle - 
\langle \phi, \mathcal{A} \omega \rangle = 
\langle \mathcal{C}_A \phi, \mathcal{B}_0\omega \rangle, 
\end{align}
and for all $\phi \in \mathbf{D}(\mathcal{D}_0^*)$ and $\omega \in \mathbf{D}(\mathcal{D})$
\begin{align}
\label{eq:CD}
\langle \mathcal{D}_0^* \phi, \omega \rangle - 
\langle \phi, \mathcal{D} \omega \rangle = 
\langle \mathcal{C}_D\phi, \mathcal{B}_0\omega \rangle.
\end{align}

A weak solution of \eqref{eq:bc_sys} then satisfies for all $\phi \in \mathbf{D}(\mathcal{A}_0^*)$
\begin{multline}
\langle \phi, \ddot{\omega}(t) \rangle + \langle \mathcal{D}^*_0\phi, \dot{\omega}(t) \rangle + \langle \mathcal{A}_0^* \phi, \omega(t) \rangle = \\
\langle \mathcal{C}_A \phi, u(t) - \mathcal{B}_1 \dot{\omega}(t) \rangle + \langle \mathcal{C}_D \phi, \dot{u}(t) - \mathcal{B}_1 \ddot{\omega}(t) \rangle, \label{eq:bc_sys-weak} 
\end{multline}
Note that the damping operator $\mathcal{D}$ may introduce a term proportional to $\dot{u}(t)$ through boundary terms captured by the operator $\mathcal{C}_D$ defined in \eqref{eq:24}. Such a term was also considered in \eqref{eqm-1} where it lead to a jump in the velocity of the step response, see \eqref{eqm-2}. The formulation (\ref{eq:bc_sys-weak}) should also be interpreted in the weak sense with respect to time. This weak form is the correct description of the system.

For the step response, i.e, the solution for $u(t) = u_0$ for $t > 0$ and $u(t)=0$ for negative time, it is assumed that the system was at rest till time zero. Since the system is causal, $\omega(t) = 0$ for $t < 0$, i.e., the initial conditions in (\ref{eq:bc_sys}) are regarded as values at $t=0^-$. Taking $u(t)$ equal to the step function in (\ref{eq:bc_sys-weak}) and integrating from $t=0^-$ to $t$ yields
\begin{multline}
\label{eq:bc_sys_wi}
\langle \phi, \dot{\omega}(t) \rangle + \langle \mathcal{D}^*_0\phi, \omega(t) \rangle + \langle \mathcal{A}_0^* \phi, W(t) \rangle = \\
\langle \mathcal{C}_A \phi, u_0t - \mathcal{B}_1\omega(t) \rangle + \langle \mathcal{C}_D \phi, u_0 - \mathcal{B}_1 \dot{\omega}(t) \rangle, 
\end{multline}
where $W(t) = \int_{0^-}^t \omega(\tau) d\tau$. Letting $t$ approach zero from above and assuming that $\omega(t)$ is continuous at $t = 0$, it follows that $\langle \phi, \dot{\omega}(0^+) \rangle = \langle  \mathcal{C}_D\phi, u_0- \mathcal{B}_1\dot{\omega}(0^+)\rangle$ for all $\phi \in \mathbf{D}(\mathcal{A}_0^*)$. So this gives a relation on $\dot{\omega}(0^+)$. 

After integrating (\ref{eq:bc_sys_wi}) once more, all derivatives on the left hand side  are gone, and since $\omega(t)$ and thus its antiderivatives take values in  $X$, it is well-defined. However, on the right hand side the term $\langle \mathcal{C}_D \phi, u_0t - \mathcal{B}_1 \omega(t) \rangle$ remains, and so even for a weak solution, this must have a meaning. In general this will imply that $\omega(t) \in {\mathbf D}({\mathcal B}_1)$ or that $\mathcal{C}_D = 0$.

Now it is not hard to see that our step response will be the weak solution on $(0^+,\infty)$ of 
\begin{align}
\langle \phi, & \ddot{\omega}(t) \rangle + \langle \mathcal{D}^*_0\phi, \dot{\omega}(t) \rangle + \langle \mathcal{A}_0^* \phi, \omega(t) \rangle = \nonumber \\
&\quad 
\langle \mathcal{C}_A \phi, u_0 - \mathcal{B}_1 \dot{\omega}(t) \rangle - \langle \mathcal{C}_D \phi, \mathcal{B}_1 \ddot{\omega}(t) \rangle,  
\label{eq:bc_sys-weak2} 
\\
&\omega(0) = 0, \quad  \langle \phi, \dot{\omega}(0^+) \rangle = \langle  \mathcal{C}_D\phi, u_0- \mathcal{B}_1\dot{\omega}(0^+)\rangle. 
\nonumber
\end{align}
From this it follows that $ \dot{\omega}(0^+)$ may be unequal to zero when $\mathcal{C}_D \neq 0$, even when $ \mathcal{B}_1=0$. The initial condition on $\dot{\omega}(0^+)$ could be rewritten in a strong from by using adjoints, but since ${\mathcal C}_D$ can be a point evaluation, this does not bring much. In order to keep the notation simple, (\ref{eq:bc_sys-weak2}) will be used in the following, but this equation needs to be understood in its integrated from, see (\ref{eq:bc_sys_wi}).

\begin{theorem} 
\label{thm:boundary_control}
 Let ${\mathcal A}_0$, as defined in (\ref{eq:23}), be a closed, densely defined operator, and let $u_0 \in U$ be given. 
A weak solution $\omega(t)$ of \eqref{eq:bc_sys-weak2} of the form \eqref{eq:omega} exists if and only if 
\begin{equation}
\label{eq:25a}
  {\mathcal A}_0 \omega_{2} = 0,
\end{equation}
for all $\phi \in \mathbf{D}(\mathcal{A}_0^*)$ there holds
\begin{align}
   \label{eq:25b-w}
   \langle {\mathcal D}_0^*\phi, \omega_{2}\rangle  + \langle {\mathcal A}_0^*\phi, \omega_{1}\rangle =&\ -\langle {\mathcal C}_A \phi, \mathcal{B}_1 \omega_2 \rangle, \\
  \langle \phi, \omega_{2}\rangle + \langle {\mathcal D}_0^*\phi, \omega_{1}\rangle  +\langle  {\mathcal A}_0^*\phi, \omega_{0}\rangle =&\  \nonumber \\
 \langle {\mathcal C}_A \phi, u_0 - {\mathcal B}_1  \omega_1\rangle - \langle {\mathcal C}_D \phi,  {\mathcal B}_1 \omega_2\rangle, \label{eq:25c-w}
\end{align}
and 
\begin{align}
  \langle \phi, \ddot{\omega}_{st}(t) \rangle +  \langle&  \mathcal{D}^*_0\phi,  \dot{\omega}_{st}(t) \rangle + \langle \mathcal{A}_0^* \phi, \omega_{st}(t) \rangle =   \nonumber \\
  -& \langle \mathcal{C}_A \phi,  \mathcal{B}_1 \dot{\omega}_{st}(t) \rangle - \langle \mathcal{C}_D \phi, \mathcal{B}_1 \ddot{\omega}_{st}(t) \rangle,  \label{eq:25d}
\\
\nonumber
\omega_{st}(0)=&\ -\omega_{0}, \\
\nonumber \langle \phi,\omega_1+ \dot{\omega}_{st}(0^+) \rangle =&\  \langle  \mathcal{C}_D\phi, u_0- \mathcal{B}_1(\omega_1 +\dot{\omega}_{st}(0^+))\rangle. 
%
%  0 =&\ \ddot{\omega}_{st}(t) + {\mathcal D} \dot{\omega}_{st}(t) + {\mathcal A} \omega_{st}(t), \quad {\mathcal B}_0\omega_{st}(t) + \mathcal{B}_1 \dot{\omega}_{st}(t) =0, \\
%  & \omega_{st}(0)=-\omega_{0}, \quad  \dot{\omega}_{st}(0)= (I + \mathcal{C}_D^*\mathcal{B}_1)^{-1}\mathcal{C}^*_Du_0 -\omega_{1}. 
\end{align}
with $\omega_{st}(t) \rightarrow 0$ for $t \rightarrow \infty$. 
 In particular, $\omega_2 \in \mathbf{D}(\mathcal{A}_0)$. 

Furthermore, the Laplace transform $\Omega(s)$ of $\omega(t)$ equals $G(s)u_0/s$, where $G(s)$ is the transfer function of 
(\ref{eq:bc_sys-weak}). Hence, if for every $u_0 \in U$ such a solution exists, then the transfer function can be written as in \eqref{eq:Gs} and \eqref{eq:Gk} holds.
\end{theorem}
{\bf Proof}  
Substituting (\ref{eq:omega}) into (\ref{eq:bc_sys-weak}) and equating the $t^k$-terms, it follows that for all $\phi \in \mathbf{D}(\mathcal{A}_0^*)$
\begin{align}
\label{eq:25a-w}
  \langle {\mathcal A}_0^*\phi, \omega_{2}\rangle =&\ 0\\
  \label{eq:25b-w-bis}
   \langle {\mathcal D}_0^*\phi, \omega_{2}\rangle  + \langle {\mathcal A}_0^*\phi, \omega_{1}\rangle =&\ -\langle {\mathcal C}_A \phi, \mathcal{B}_1 \omega_2 \rangle, \\
   %\nonumber 
  \langle \phi, \omega_{2}\rangle + \langle {\mathcal D}_0^*\phi, \omega_{1}\rangle  +\langle  {\mathcal A}_0^*\phi, \omega_{0}\rangle =&\ %\\
  \langle {\mathcal C}_A \phi, u_0 - {\mathcal B}_1 \omega_1\rangle - \langle {\mathcal C}_D \phi,  {\mathcal B}_1 \omega_2\rangle,  \label{eq:25c-w-bis}
\end{align}
and (\ref{eq:25d}). 
Since (\ref{eq:25a-w}) holds for all $\phi \in {\mathbf D}({\mathcal A}_0^*)$, and since ${\mathcal A}_0$ is closed, $\omega_2 \in {\mathbf D}({\mathcal A}_0)$ and (\ref{eq:25a}) holds. 

For the second part, note that the relation for the transfer function of the system given by (\ref{eq:bc_sys-weak}). Using the definition of the transfer function $G(s)$, see \cite[Section 10.2]{Ruth}, $G(s)u_0$ satisfies
\begin{align}
\nonumber 
\langle \left(s^2 + s \mathcal{D}^*_0 + \mathcal{A}^*_0 \right) \phi, G(s)u_0 \rangle &= \langle \mathcal{C}_A \phi, u_0 - s\mathcal{B}_1G(s)u_0 \rangle \\
 + \langle \mathcal{C}_D & \phi, su_0- s^2\mathcal{B}_1G(s)u_0  \rangle. 
\label{eq:bc_TF}
\end{align}
Taking the Laplace transform of the proposed solution $\omega(t)$ in \eqref{eq:omega} yields
\begin{equation*}
  \Omega(s) = \frac{\omega_2}{s^3} + \frac{\omega_1}{s^2} + \frac{\omega_0}{s} + \Omega_{st}(s),
\end{equation*}
where $\Omega_{st}(s)$ is the Laplace transform of $\omega_{st}(t)$. 

Similarly as in Theorem \ref{Th-1}, it suffices to show that $\Omega(s)$ equals $G(s)u_0/s$, or equivalently, $s\Omega(s)$ satisfies the equation for $G(s)u_0$ above. 
Note that
\begin{align}
\nonumber
&\langle \left(s^2 + s \mathcal{D}^*_0 + \mathcal{A}^*_0 \right) \phi,  \Omega(s) \rangle \\ \nonumber
&= \left\langle \left(s^2 + s \mathcal{D}^*_0 + \mathcal{A}^*_0 \right) \phi,  \frac{\omega_2}{s^3} + \frac{\omega_1}{s^2} + \frac{\omega_0}{s} + \Omega_{st}(s) \right\rangle 
\end{align}
\begin{align}
\nonumber
&= \frac{\langle \mathcal{A}^*_0 \phi, \omega_2 \rangle}{s^3} + 
\left[ \langle \mathcal{D}^*_0\phi, \omega_2 \rangle + \langle \mathcal{A}^*_0 \phi, \omega_1 \rangle \right] \frac{1}{s^2} \\ \nonumber 
& \qquad\qquad  + \left[ \langle \phi, \omega_2 \rangle + \langle \mathcal{D}^*_0\phi, \omega_1 \rangle + \langle \mathcal{A}^*_0 \phi, \omega_0 \rangle \right] \frac{1}{s} \\
\nonumber
& \qquad\qquad + \langle \phi, \omega_1 \rangle + \langle \mathcal{D}^*_0 \phi, \omega_0 \rangle + s \langle \phi, \omega_0 \rangle \\ \nonumber 
& \qquad\qquad + \langle \left(s^2 + s \mathcal{D}_0^* + \mathcal{A}^*_0 \right) \phi, \Omega_{st}(s) \rangle \\
\nonumber
&= \frac{-\langle \mathcal{C}_A \phi,  \mathcal{B}_1 \omega_2 \rangle}{s^2} + \frac{-\langle \mathcal{C}_D \phi, \mathcal{B}_1 \omega_2 \rangle + \langle \mathcal{C}_A \phi, u_0 - \mathcal{B}_1 \omega_1 \rangle}{s} \\ 
&\qquad\qquad + \langle \phi, \omega_1 \rangle + \langle \mathcal{D}^*_0 \phi, \omega_0 \rangle + s \langle \phi, \omega_0 \rangle \nonumber \\
&\qquad\qquad + \langle \left(s^2 + s \mathcal{D}_0^* + \mathcal{A}^*_0 \right) \phi, \Omega_{st}(s) \rangle, \label{eq:38}
\end{align}
where \eqref{eq:25a-w}, \eqref{eq:25b-w} and \eqref{eq:25c-w} have been used. 

Similar, as in the proof of Theorem \ref{Th-1} we define the Lapace transform as $\int_{0^+}^{\infty} f(t) e^{-st} dt$. Taking the Laplace transform of (\ref{eq:25d}) gives
\begin{align}
  & \langle \phi, s^2 \Omega_{st}(s) - s \omega_{st}(0) - \dot{\omega}_{st}(0^+)\rangle \nonumber \\
   &\qquad\quad + \langle  \mathcal{D}^*_0\phi, s \Omega_{st}(s)- \omega_{st}(0)\rangle 
   +
   \langle \mathcal{A}^*_0 \phi, \Omega_{st}(t) \rangle \nonumber \\
   &= - \langle \mathcal{C}_A \phi, \mathcal{B}_1 [s \Omega_{st}(s) - \omega_{st}(0) ]\rangle \nonumber \\
   & \qquad\quad   - \langle \mathcal{C}_D \phi, \mathcal{B}_1 [s^2 \Omega_{st}(s) - s \omega_{st}(0) - \dot{\omega}_{st}(0^+)]\rangle. \label{eq:39}
\end{align}
This leads to the following equality for $\Omega_{st}$
\begin{align*}
\langle   (s^2& +  s\mathcal{D}_0^* + \mathcal{A}_0^* )\phi, \Omega_{st}(s) \rangle  \\
=&\ 
\langle \phi, s\omega_{st}(0) + \dot{\omega}_{st}(0^+) \rangle 
+ \langle \mathcal{D}_0^* \phi, \omega_{st}(0) \rangle \\
&\qquad\qquad - \langle \mathcal{C}_A \phi,  \mathcal{B}_1\left(s\Omega_{st}(s) - \omega_{st}(0) \right)\rangle\\
&\qquad\qquad - \langle \mathcal{C}_D \phi,  \mathcal{B}_1\left(s^2\Omega_{st}(s) - s\omega_{st}(0) - \dot{\omega}_{st}(0^+) \right)\rangle 
 \\
=& -s\langle \phi, \omega_0 \rangle - \langle \phi, \omega_1 \rangle + \langle \mathcal{C}_D \phi, u_0 - {\mathcal B}_1 \left(\omega_1 +\dot{\omega}_{st}(0^+)\right) \rangle  \\
&\qquad\qquad - \langle \mathcal{D}_0^* \phi, \omega_0 \rangle - \langle \mathcal{C}_A \phi,  \mathcal{B}_1\left(s\Omega_{st}(s) + \omega_0 \right)\rangle\\
&\qquad\qquad - \langle \mathcal{C}_D \phi,  \mathcal{B}_1\left(s^2\Omega_{st}(s) + s\omega_0 - \dot{\omega}_{st}(0^+) \right)\rangle  \\
=& -s\langle \phi, \omega_0 \rangle - \langle \phi, \omega_1 \rangle - \langle \mathcal{D}_0^* \phi, \omega_0 \rangle \\
&\qquad\qquad +\langle \mathcal{C}_D \phi, u_0 - \mathcal{B}_1 (\omega_1 +  s\omega_0 + s^2 \Omega_{st}(s)) \rangle \\
&\qquad\qquad - \langle \mathcal{C}_A \phi,  \mathcal{B}_1\left(s\Omega_{st}(s) + \omega_0 \right)\rangle,
\end{align*}
where the initial conditions in \eqref{eq:25d} have been used. Substituting this into equation (\ref{eq:38}) gives
\begin{align*}
\langle (s^2 &  
+ s \mathcal{D}^*_0 + \mathcal{A}^*_0 ) \phi,  \Omega(s) \rangle \\
=&
 \frac{-\langle \mathcal{C}_A \phi,  \mathcal{B}_1 \omega_2 \rangle}{s^2} + \frac{-\langle \mathcal{C}_D \phi, \mathcal{B}_1 \omega_2 \rangle + \langle \mathcal{C}_A \phi, u_0 - \mathcal{B}_1 \omega_1 \rangle}{s} \\
 &\qquad\quad + \langle \mathcal{C}_D \phi, u_0 - \mathcal{B}_1 (\omega_1 +  s\omega_0 + s^2 \Omega_{st}(s)) \rangle \\
&\qquad\quad - \langle \mathcal{C}_A \phi,  \mathcal{B}_1\left(s\Omega_{st}(s) + \omega_0 \right)\rangle\\
=& \left\langle \mathcal{C}_A \phi, \frac{u_0}{s} - \mathcal{B}_1 \left( \frac{\omega_2}{s^2} + \frac{\omega_1}{s} + \omega_0 + s\Omega_{st}(s) \right) \right\rangle \\
&\qquad\quad +  \langle \mathcal{C}_D\phi, u_0- \mathcal{B}_1 \left(  \omega_2 \tfrac{1}{s} + \omega_1 + s\omega_0 + s^2 \Omega_{st}(s) \right) \rangle \\
=&\ \left\langle \mathcal{C}_A \phi, \frac{u_0}{s} - s \mathcal{B}_1\Omega(s) \right\rangle +  \langle \mathcal{C}_D\phi, u_0- s^2 \mathcal{B}_1  \Omega(s) \rangle.
\end{align*}
Comparing this with (\ref{eq:38}), it follows that $s\Omega(s)$ satisfies the same equation as $G(s)u_0$. Since the transfer function is unique, they must be equal.
\hfill$\Box$
\medskip

The above theorem only gives weak equations for $\omega_1$ and $\omega_2$. However, by imposing an extra condition, more direct equations can be obtained.
\begin{lemma} 
\label{La:1}
  If ${\mathcal B}_1=0$, then under the conditions as stated in Theorem \ref{thm:boundary_control} there holds that $\omega_1 \in {\mathbf D}({\mathcal A}_0)$ and 
   \begin{equation}
  \label{eq:25b0}
  {\mathcal D}_0\omega_{2} +  {\mathcal A}_0\omega_{1} = 0.
   \end{equation}
\end{lemma}
{\bf Proof}  Since ${\mathcal B}_1=0$, equation (\ref{eq:25b-w}) becomes
\[
    \langle {\mathcal D}_0^*\phi, \omega_{2}\rangle  + \langle {\mathcal A}_0^*\phi, \omega_{1}\rangle =0.
\]
Since $\omega_2 \in {\mathbf D}({\mathcal A}_0) \subseteq  {\mathbf D}({\mathcal D}_0) $, this implies that 
\[
  \langle {\mathcal A}_0^*\phi, \omega_{1}\rangle = - \langle\phi, {\mathcal D}_0 \omega_{2}\rangle \quad \mbox{ for all } \phi \in {\mathbf D}({\mathcal A}_0^*).
\]
Since ${\mathcal A}_0$ is closed, this gives that $\omega_1 \in {\mathbf D}({\mathcal A}_0) $ and $\langle\phi,  {\mathcal A}_0 \omega_{1}\rangle = - \langle\phi, {\mathcal D}_0 \omega_{2}\rangle$. From this and the fact that ${\mathbf D}({\mathcal A}_0^*)$ is dense in $X$, (\ref{eq:25b}) follows.
\hfill$\Box$
\medskip

In the following lemma, another condition under which (\ref{eq:25b-w}) and (\ref{eq:25c-w}) can be simplified is given.
\begin{lemma} 
\label{La:2}
  Assume that ${\mathcal B}_0$ is surjective from ${\mathbf D}({\mathcal A})$ to $U$. Then under the conditions as stated in Theorem \ref{thm:boundary_control} there holds that $\omega_0,\omega_1 \in {\mathbf D}({\mathcal A})$, and
   \begin{align}
   \label{eq:25b}
   {\mathcal D}_0\omega_{2} +  {\mathcal A}\omega_{1}  =&\ 0, \quad {\mathcal B}_0 \omega_1 + {\mathcal B}_1 \omega_2=0,\\
  \label{eq:25c}
 \omega_{2}+ {\mathcal D} \omega_{1} +  {\mathcal A} \omega_{0} =&\ 0, \quad {\mathcal B}_0 \omega_0 + {\mathcal B}_1 \omega_1=u_0.
\end{align}
\end{lemma}
 {\bf Proof} 
 By assumption, there exists a $B_1\in \mathbf{D}(\mathcal{A})$ such that ${\mathcal B}_0 B_1 = -  \mathcal{B}_1 \omega_2$. So with (\ref{eq:CA}), (\ref{eq:25b-w}) can be rewritten as
\[
  \langle {\mathcal D}_0^*\phi, \omega_{2}\rangle  + \langle {\mathcal A}_0^*\phi, \omega_{1}\rangle =  \langle {\mathcal A}_0^*\phi, B_{1}\rangle - \langle \phi, {\mathcal A} B_{1}\rangle.
\]
Furthermore, (\ref{eq:CD}) implies that 
\[
 \langle {\mathcal D}_0^*\phi, \omega_{2}\rangle = \langle \phi, {\mathcal D} \omega_{2}\rangle + \langle C_D \phi, {\mathcal B}_0\omega_2\rangle = \langle \phi, {\mathcal D}_0 \omega_{2}\rangle,
\]
where it was used that $\omega_2 \in {\mathbf D}({\mathcal A}_0)\subseteq {\mathbf D}({\mathcal D}_0)$. Combining these gives
\[
 \langle \phi,  {\mathcal D}_0 \omega_{2}\rangle + \langle {\mathcal A}_0^* \phi, \omega_{1} - B_1\rangle = - \langle \phi, {\mathcal A}B_1 \rangle.
\]
Because ${\mathcal A}_0$ is closed, it follows that $\omega_{1} - B_1 \in {\mathbf D}({\mathcal A}_0)$, and thus that $0={\mathcal B}_0( \omega_{1} - B_{1}) = {\mathcal B}_0 \omega_1 - {\mathcal B}_0 B_{1} = {\mathcal B}_0 \omega_1 + {\mathcal B}_1\omega_2$. Furthermore, 
\[
  - {\mathcal A} B_{1}=  {\mathcal D}_0\omega_{2}  +  {\mathcal A}_0 ( \omega_{1} - B_{1}) =  {\mathcal D}_0\omega_{2}  +  {\mathcal A} \omega_{1} -  {\mathcal A}B_{1},
\]
which implies (\ref{eq:25b}).  

Since $\omega_1 \in {\mathbf D}({\mathcal A})\subseteq {\mathbf D}({\mathcal D})$, (\ref{eq:25c-w}) can be rewritten as
\begin{align*}
\langle \phi, \omega_{2}\rangle + \langle \phi, & {\mathcal D}\omega_{1}\rangle  + \langle {\mathcal C}_D \phi,  {\mathcal B}_0 \omega_1\rangle+\langle  {\mathcal A}_0^*\phi, \omega_{0}\rangle \\
&= \langle {\mathcal C}_A \phi, u_0 - {\mathcal B}_1 \omega_1\rangle - \langle {\mathcal C}_D \phi,  {\mathcal B}_1 \omega_2\rangle,
\end{align*}
which, because of  ${\mathcal B}_0 \omega_1 + {\mathcal B}_1\omega_2=0$, simplifies to 
\[
  \langle \phi, \omega_{2}\rangle + \langle \phi, {\mathcal D}\omega_{1}\rangle  +\langle  {\mathcal A}_0^*\phi, \omega_{0}\rangle = \langle {\mathcal C}_A \phi, u_0 - {\mathcal B}_1 \omega_1\rangle .
\]
Equation (\ref{eq:25c}) is obtained in a similar way. In particular, choose $B_{0}\in {\mathbf D}({\mathcal A})$ such that ${\mathcal B}_0 B_{0}=u_0- {\mathcal B}_1 \omega_1$, and write 
\[
  \langle {\mathcal C}_A \phi, u_0 - {\mathcal B}_1 \omega_1\rangle = \langle  {\mathcal A}_0^*\phi, B_{0}\rangle - \langle  \phi, {\mathcal A} B_{0}\rangle.
\]
This shows that $\omega_0 - B_0 \in {\mathbf D}({\mathcal A}_0)$ which implies that ${\mathcal B}_0 \omega_0 = u_0- {\mathcal B}_1 \omega_1$ and (\ref{eq:25c}).
\hfill$\Box$
\medskip

Similar as in Theorem \ref{thm:om1}, more conditions on $\omega_2$, $\omega_1$, and $\omega_0$ can be obtained when ${\mathcal A}_0$ and ${\mathcal D}_0$ are symmetric. Two cases are distinguished: ${\mathcal B}_1 =0$ and $\mathcal{C}_D = 0$.
\begin{theorem} 
\label{thm:om1-bc}
Assume that ${\mathcal A}_0$, ${\mathcal D}_0$ are symmetric, and ${\mathcal D}_0$ is non-negative. Then the following holds
\begin{itemize}
\item When ${\mathcal B}_1=0$, then 
  \begin{itemize}
  \item ${\mathcal D}_0\omega_2=0$ and ${\mathcal A}_0\omega_1=0$;
  \item for every $\phi \in \ker(\mathcal{A}_0) \cap \ker(\mathcal{D}_0) $ there holds
\begin{align}
\label{eq:thm2-bc1}
\langle \phi, \omega_2 \rangle =&\ \langle \mathcal{C}_A \phi, u_0  \rangle;\\
\label{eq:thm2-bc2}
\langle \phi, \omega_1 \rangle =&\ \langle \mathcal{C}_D \phi, u_0  \rangle;\\
\label{eq:thm2-bc3}
\langle \phi, \omega_0 \rangle =&\ 0.
%\label{eq:thm2-bc1}
\end{align}
  \end{itemize}
  
  \item When ${\mathcal C}_D=0$, then 
 \begin{itemize}
  \item for every $\phi \in \ker(\mathcal{A}_0) \cap \ker(\mathcal{D}_0)$ there holds
\begin{align}
\label{eq:thm2-bc1-a}
\langle \phi, \omega_2 \rangle =&\ \langle \mathcal{C}_A \phi, u_0 -{\mathcal B}_1\omega_1 \rangle;\\
\label{eq:thm2-bc2-a}
\langle \phi, \omega_1 \rangle =&\ - \langle \mathcal{C}_A \phi, {\mathcal B}_1\omega_0  \rangle.
\end{align}
\item
 for every $\phi \in \ker(\mathcal{A}_0) \cap \ker(\mathcal{D}_0) \cap \ker({\mathcal C}_A)$
\begin{equation}
\label{eq:thm2-bc3-a}
\langle \phi, \omega_0 \rangle = 0.
\end{equation}
  \end{itemize}
\end{itemize}
\end{theorem}
{\bf Proof} 
Because $\mathcal{A}_0$ and $\mathcal{D}_0$ are symmetric, $\mathbf{D}(\mathcal{A}_0) \subseteq \mathbf{D}(\mathcal{A}_0^*)$ and $\mathbf{D}(\mathcal{D}_0) \subseteq \mathbf{D}(\mathcal{D}_0^*)$. 

\noindent
{\bf Case 1:} $\mathcal{B}_1=0$. Because $\omega_2 \in \mathbf{D}(\mathcal{A}_0)\subseteq \mathbf{D}(\mathcal{A}_0^*)$, \eqref{eq:25b-w} can be evaluated with $\phi=\omega_2$. This shows that
\begin{align*}
0 &= \langle {\mathcal D}_0^* \omega_2, \omega_2 \rangle + \langle {\mathcal A}_0^* \omega_2, \omega_1\rangle \\
 &=  \langle {\mathcal D}_0 \omega_2,  \omega_2 \rangle + \langle {\mathcal A}_0 \omega_2, \omega_1\rangle \\
&= \langle {\mathcal D}_0 \omega_2,  \omega_2 \rangle 
\end{align*}
where the first equality holds by the symmetry and last equality follows from (\ref{eq:25a}).
Since  $\mathcal{D}_0$ is non-negative by assumption, ${\mathcal D}_0\omega_2=0$. Now (\ref{eq:25b0}) gives that ${\mathcal A}_0\omega_1=0$.

Taking $\phi \in \ker(\mathcal{A}_0) \cap \ker(\mathcal{D}_0)$ in \eqref{eq:25c-w}, \eqref{eq:thm2-bc1} follows.

For $\phi \in \ker(\mathcal{A}_0) \cap \ker(\mathcal{D}_0)$, define
\[
  v_{st}(t) = \langle \phi, \omega_{st}(t) \rangle.
\]
Equation \eqref{eq:25d} then shows that $\ddot{v}_{st}(t) =0$. Since $\omega_{st}(t)$ is stable, $v_{st}(0)=\dot{v}_{st}(0^+) = 0$. The condition in \eqref{eq:thm2-bc3} now follows from the requirement that $v_{st}(0) = 0$ and the initial condition in \eqref{eq:25d}. Similarly, the second initial condition in \eqref{eq:25d} combined with $\langle \phi, \dot{\omega}_{st}(0^+) \rangle = \dot{v}_{st}(0^+) = 0$ gives \eqref{eq:thm2-bc2}. 
\smallskip

\noindent
{\bf Case 2: $\mathcal{C}_D = 0$}. 
Taking $\phi \in \ker(\mathcal{A}_0) \cap \ker(\mathcal{D}_0)$ in \eqref{eq:25c-w}, \eqref{eq:thm2-bc1-a} follows. 

For \eqref{eq:thm2-bc2-a}, note that (\ref{eq:39}) for this $\phi$ becomes
\[
  \langle \phi, s^2 \Omega_{st}(s) - s \omega_{st}(0) - \dot{\omega}_{st}(0^+) \rangle %\\
  = - \langle {\mathcal C}_A \phi, {\mathcal B}_1(s \Omega_{st}(s) - \omega_{st}(0) ) \rangle.
\]
Applying the boundary conditions of (\ref{eq:25d}) gives 
\begin{equation}
\label{eq:49}
  \langle \phi, s^2 \Omega_{st}(s) + s \omega_0 \rangle  + \langle {\mathcal C}_A \phi, {\mathcal B}_1s \Omega_{st}(s)\rangle %\\
  %\nonumber
  = \langle \phi, \omega_1 \rangle  - \langle {\mathcal C}_A \phi, {\mathcal B}_1 \omega_0 \rangle.
\end{equation}
Since $\Omega_{st}$ is stable, the value at $s=0$ exists (or at least $\lim_{s \downarrow 0}$ exists). Taking this limit in the above equality gives (\ref{eq:thm2-bc2-a}). In particular, we know that the sum in (\ref{eq:49}) equals zero. So for $\phi \in \ker(\mathcal{A}_0) \cap \ker(\mathcal{D}_0) \cap \ker( {\mathcal C}_A)$ this equality implies that
\[
  \langle \phi, s^2 \Omega_{st}(s) + s \omega_0 \rangle  =0.
\]
Dividing by $s$, taking the limit as $s\downarrow 0$, and using once more that $\lim_{s \downarrow 0}\Omega_{st}(s)$ exists, the above equality gives \eqref{eq:thm2-bc3-a}.
\hfill$\Box$
\medskip

Similar as before, an explicit expression for $\omega_2$, $\omega_1$, and $\omega_0$ can be obtained when $\mathcal{A}_0$ and $\mathcal{D}_0$ are diagonalizable on the same basis of eigenfunctions as in \eqref{eq:A0D0_spectral}, 
%$u_0$ lies in the range of $\mathcal{B}_0 : \mathbf{D}(\mathcal{B}_0) \rightarrow U$, 
and $\mathcal{B}_1 = 0$. 
\begin{theorem} 
\label{thrm1-bc}
%\normalfont  
Let $\mathcal{B}_1 = 0$ and let ${\mathcal A}_0$ and $\mathcal{D}_0$ be non-negative self-adjoint operators with compact resolvent satisfying $\ker(\mathcal{D}_0) \subseteq \ker(\mathcal{A}_0)$. Let the eigenvalues and normalized eigenfunctions of $\mathcal{A}_0$ be $\lambda_n$ and $\varphi_{n}$, $n \in {\mathbb N}$. Let $\mathcal{D}_0$ have the same eigenfunctions $\varphi_n$ and let its eigenvalues be $\mu_n$, $n \in \mathbb{N}$. 
%If there exists a $b_0 \in \mathbf{D}(\mathcal{B}_0)$ such that $\mathcal{B}_0 b_0 = u_0$, then the 
Then the $\omega_0$, $\omega_1$, and $\omega_2$ satisfying \eqref{eq:25a}--(\ref{eq:25d}) are
\begin{equation} 
\label{w22-bc}
   \omega_{2}=\sum_{n\in \mathbb{J}}\alpha_n \varphi_n, \quad 
   \omega_1= \sum_{n\in {\mathbb I}}\beta_n \varphi_n, \quad
   \omega_{0}= \sum_{n\in \mathbb{N} \backslash \mathbb{J}} \gamma_n \varphi_n,
\end{equation}
where $\mathbb{I} = \{ n \mid \lambda_n = 0 \}$, $\mathbb{J} = \{ n \mid \mu_n = 0 \}$, and
\begin{align}
\label{alfa+beta-bc}
        \alpha_{n} &=\langle \mathcal{C}_A\varphi_{n}, u_{0}\rangle, \quad  
        \beta_{n}= \left\{ \begin{array}{ll}
         \langle\mathcal{C}_D \varphi_{n},u_{0}\rangle & n \in \mathbb{J} \\
        \frac{\langle \mathcal{C}_A\varphi_{n},u_{0}\rangle}{\mu_n}, & n \in \mathbb{I} \backslash \mathbb{J}
        \end{array}
        \right. 
        \\
        \gamma_{n} &= %-\langle \varphi_n , b_0 \rangle + 
        \left\{ \begin{array}{ll}
        \frac{\langle \mathcal{C}_D \varphi_n, u_0 \rangle}{\mu_n} - \frac{\langle \mathcal{C}_A \varphi_n, u_0 \rangle}{\mu_n^2} & n \in \mathbb{I} \backslash \mathbb{J} \\
        \frac{\langle \mathcal{C}_A \varphi_{n}, u_{0} \rangle}{\lambda_n} & n \in \mathbb{N} \backslash \mathbb{I}
        \end{array}
        \right. . \label{gamma-bc}
\end{align}
\end{theorem}
{\bf Proof}  
Observe that $\ker(\mathcal{D}_0) \subseteq \ker(\mathcal{A}_0)$ implies that $\mathbb{J} \subseteq \mathbb{I}$. Theorem \ref{thm:om1-bc} shows that $\omega_2 \in \ker(\mathcal{D}_0)$ and $\omega_1 \in \ker({\mathcal A}_0)$. This implies $\omega_2$ is of the form in \eqref{w22-bc} for certain coefficients $\alpha_n$, and $\omega_1$ is of the form in \eqref{w22-bc} for certain coefficients $\beta_n$. Furthermore, taking $\varphi_n$ with $n \in \mathbb{J}$ in \eqref{eq:thm2-bc3} shows that $\omega_0$ is of the form \eqref{w22-bc} for certain coefficients $\gamma_n$. 

%Furthermore, $\mathcal{B}_0(\omega_0 - b_0) = 0$, which means that $\omega_0$ is of the form in \eqref{w22-bc} for certain coefficients $\gamma_n \in \mathbb{R}$. \textcolor{blue}{This is the most tricky part. I think we do not need here that $b_0 \in \mathbf{D}(\mathcal{A})$, right? }

Inserting these expressions into \eqref{eq:thm2-bc1} and taking $\phi = \varphi_n$ for $n \in \mathbb{J}$ yields the expression for $\alpha_n$ in \eqref{alfa+beta-bc}. 
The expression for $\beta_n$ with $n \in \mathbb{J}$ follows by taking $\phi = \varphi_n$ in \eqref{eq:thm2-bc2}.  

Taking $\phi = \varphi_n$ for $n \in \mathbb{I} \backslash \mathbb{J}$ in \eqref{eq:25c-w} yields the expression for $\beta_n$ in \eqref{alfa+beta-bc} for $n \in \mathbb{I} \backslash \mathbb{J}$. 

Using (\ref{eq:25c-w}) once more, but now with $\phi= \varphi_n$ for $n \not\in \mathbb{I}$, gives the coefficients $\gamma_n$ of (\ref{gamma-bc}) for these $n$'s. Thus it remains to find the $\gamma_n$ for $n \in {\mathbb I} \backslash {\mathbb J}$. 

%Condition \eqref{eq:thm2-bc3} in Theorem \ref{thm:om1-bc} with $\phi = \varphi_n$ for $n \in \mathbb{J}$ yields the expression for $\gamma_n$ with $n \in \mathbb{J}$ in \eqref{gamma-bc}. Taking $\phi = \varphi_n$ for $n \in \mathbb{N} \backslash \mathbb{I}$ in \eqref{eq:25c-w} yields the expression for $\gamma_n$ with $n \in \mathbb{N} \backslash \mathbb{I}$ in \eqref{gamma-bc}. 

To determine the expressions for $\gamma_n$ with $n \in \mathbb{I} \backslash \mathbb{J}$, consider $v_n(t) = \langle \varphi_n, \omega_{st}(t) \rangle$. From the differential equation \eqref{eq:25d} and the expressions (\ref{w22-bc}), it follows that
\begin{align}
  & \ddot{v}_n(t) + \mu_n \dot{v}_n(t) = 0, \\ 
  & v_n(0) = - \langle \varphi_n, \omega_0 \rangle, \quad \dot{v}_n(0^+) = \langle \mathcal{C}_D \varphi_n, u_0 \rangle - \langle \varphi_n, \omega_1 \rangle. \nonumber 
\end{align}
The solution to this ODE is given by
\begin{equation*}
v_n(t) = v_n(0) + \frac{\dot{v}_n(0^+)}{\mu_n}(1 - e^{-\mu_n t} ). 
\end{equation*}
In order to have that $\omega_{st}(t) \rightarrow 0$ for $t\rightarrow \infty$, it is thus required that
\begin{equation*}
0 = v_n(0) + \frac{\dot{v}_n(0^+)}{\mu} = -\langle \varphi_n, \omega_0 \rangle + \frac{\langle \mathcal{C}_D \varphi_n, u_0 \rangle - \langle \varphi_n, \omega_1 \rangle}{\mu_n}. 
\end{equation*}
So for $n \in \mathbb{I} \backslash \mathbb{J}$
\begin{align*}
  \gamma_n = \langle \varphi_n, \omega_0 \rangle &= \frac{\langle \mathcal{C}_D \varphi_n, u_0 \rangle - \langle \varphi_n, \omega_1 \rangle}{\mu_n} \\
  &= \frac{\langle \mathcal{C}_D \varphi_n, u_0 \rangle}{\mu_n} - \frac{\langle \mathcal{C}_A \varphi_n, u_0 \rangle}{\mu_n^2} ,
\end{align*}
where (\ref{w22-bc}) and (\ref{alfa+beta-bc}) was used.
\hfill$\Box$
\medskip

\section{Examples} 
\label{sec:5}
In this section the time domain approach to determine the compliance is demonstrated in three examples. 

\subsection{Vibrating string}
\label{Ex:wave}
Consider the following model of a vibrating string with structural damping
    \begin{align}
       \rho(x) \frac{\partial^2 \omega}{\partial t^2}(x,t)=&\ \frac{\partial}{\partial x}\left( T(x)\frac{\partial \omega}{\partial x}(x,t) + d(x) \frac{\partial^2 \omega}{\partial t \partial x}(x,t) \right), \nonumber \\
       T(0)\frac{\partial \omega}{\partial x}(0,t)=&\ 0, \qquad
       T(L)\frac{\partial \omega}{\partial x}(L,t)=u(t), \label{eqm-3-20}
    \end{align}
where $x \in (0,L)$ and $t > 0$, $\omega(x,t)$ is the transversal displacement of the string at position $x$ and time $t$, $0 < \underline{\rho} \leq \rho(x) \leq \bar{\rho}$ is the mass density, $0 < \underline{T} \leq T(x) \leq \bar{T}$ is the tension in the string, $0 < \underline{d} \leq d(x) \leq \bar{d}$ describes the distribution of the damping, which is assumed to be smooth. The applied force at the endpoint $x = L$ is the input $u(t)$. 

Note that this model is a boundary control system of the form \eqref{eq:bc_sys} on the state space $X = L^2(0,L)$ with inner product 
\begin{equation}
  \langle f,g \rangle _X = \int_0^L f(x) \rho(x) g(x) dx.
\end{equation}
The (unbounded) operators ${\mathcal A}$ and ${\mathcal D}$ are given by
\begin{equation}
 {\mathcal A} = - \frac{1}{\rho(x)} \frac{d}{d x}\left( T(x) \frac{d}{d x} \right), \quad {\mathcal D} = - \frac{1}{\rho(x)} \frac{d}{d x}\left( d(x) \frac{d}{d x} \right), \nonumber 
\end{equation}
with domains $\mathbf{D}({\mathcal A}) = {\mathbf D}({\mathcal D}) = \{ w\in H^2(0,L)\mid \frac{dw}{dx}(0)=0\}$. 
The boundary operators ${\mathcal B}_0$ and $\mathcal{B}_1$ (with domains $\mathbf{D}(\mathcal{B}_0) = \mathbf{D}(\mathcal{A})$ and $\mathbf{D}(\mathcal{B}_1) = X$) map into the input space $U = \mathbb{R}$ and are defined as
\begin{equation}
\label{eq:59}
   {\mathcal B}_0 \omega = T(L) \frac{d \omega}{d x}(L), \qquad \mathcal{B}_1 \omega = 0.
\end{equation}
To compute the compliance, relations (\ref{eq:CA}) and (\ref{eq:CD}) need to be verified.
For $g \in {\mathbf D}({\mathcal A})$ and $f \in {\mathbf D}({\mathcal A}_0) = {\mathbf D}({\mathcal A}) \cap \ker({\mathcal B}_0)$ it is easy to show that
\begin{align}
  \langle f, {\mathcal A} g \rangle_X &= - f(L) T(L)  \frac{d g}{d x}(L) + \int_0^L  \frac{d f}{d x}(x) T(x) \frac{dg}{d x}(x)dx \nonumber \\
  &= - {\mathcal C}_A(f) {\mathcal B}_0g +  \langle {\mathcal A_0} f, g \rangle_X, 
\label{eq:62}
\end{align}
where ${\mathcal C}_A(f)=f(L)$. This shows among others that ${\mathcal A}_0$ is symmetric. 
Similarly, for $g \in {\mathbf D}({\mathcal D})=  {\mathbf D}({\mathcal A}) $ and $f \in {\mathbf D}({\mathcal D}_0) = {\mathbf D}({\mathcal A}_0)$ 
\begin{align}
  \langle f, {\mathcal D} g \rangle_X &= - f(L) d(L)  \frac{d g}{d x}(L) + \int_0^L  \frac{d f}{d x}(x) d(x) \frac{dg}{d x}(x)dx \nonumber \\
  &= - {\mathcal C}_D(f) {\mathcal B}_0g +  \langle {\mathcal D_0} f, g \rangle_X, 
   \label{eq:63}
\end{align}
where ${\mathcal C}_D(f)=\frac{d(L)}{T(L)}f(L)$. This shows among others that ${\mathcal D}_0$ is symmetric and non-negative. It is well-known that ${\mathcal A}_0$ and ${\mathcal D}_0$ are self-adjoint. 

By considering the function $x^2$ in $\mathbf{D}(\mathcal{A})$, it is easy to check that $\mathcal{B}_0 : \mathbf{D}(\mathcal{A}) \rightarrow U = \mathbb{R}$ is surjective. 

By Lemma \ref{La:2} it thus follows that $\omega_1, \omega_2 \in \mathbf{D}(\mathcal{A})$ and that $\omega_0$, $\omega_1$, and $\omega_2$ are strong solutions of \eqref{eq:25a}, \eqref{eq:25b}, and \eqref{eq:25c}. 

To determine all solution to $\mathcal{A}_0\omega_2 = 0$, note that ${\mathcal A}\omega_2 = 0$ implies that
\begin{equation} 
  \frac{d}{d x} \left(
  T(x)\frac{d \omega_2}{d x}(x)
  \right) = 0, \ \Rightarrow \ \omega_2(x) = C_1 \int_0^x \frac{1}{T(\xi)} \ \mathrm{d}\xi  + C_2, \nonumber
\end{equation}
for certain constants $C_1$ and $C_2$. 
Because $\frac{d \omega_2}{dx}(x) = C_1/T(x)$, the requirement that $\omega_2 \in \mathbf{D}(\mathcal{A})$ gives that $C_1 = 0$ and thus $\omega_2(x) = C_2$ is constant. From (\ref{eq:59}), it follows that $\omega_2(x) = C_2$ satisfies ${\mathcal B}_0 \omega_2 = 0$ for all $C_2$. 

To determine $\omega_1(x)$, note that Theorem \ref{thm:om1-bc} shows that $\mathcal{A}_0 \omega_1 = 0$. It thus follows similarly as for $\omega_2$ that $\omega_1(x) = C_3$ is constant.  

Equation (\ref{eq:25c}) now reduces to $0 = \omega_2  + {\mathcal A}\omega_{0}$ and ${\mathcal B}_0\omega_{0} = u_0$. Note that ${\mathcal A}\omega_{0}=-\omega_2$ gives that
\begin{equation}
 \label{eq:wave_om0_step1}
  %\frac{d}{dx} \left(T(x)\frac{d\omega_{0}}{d x}(x) \right) = \omega_2 \rho(x) \quad \Rightarrow
  %\quad 
  \omega_{0}(x) = \int_0^x \frac{\omega_2 \int_0^\xi \rho(\eta) \, \mathrm{d}\eta + C_4}{T(\xi)} \,\mathrm{d}\xi + C_5, \nonumber
\end{equation}
for some constants $C_4$ and $C_5$. Because $\frac{d \omega_0}{dx}(x) = (\omega_2 \int_0^x \rho(\eta) \ \mathrm{d}\eta + C_4)/T(x)$,  the boundary condition in ${\mathbf D}({\mathcal A})$ and the condition ${\mathcal B}\omega_{0} = u_0$ imply that $C_4 = 0$ and $\omega_2 = u_0 / \int_0^L \rho(\eta) \, \mathrm{d}\eta$. So $\omega_2$ has been determined.

The constants $C_5$ and $\omega_1 = C_3$ are determined from the conditions in Theorem \ref{thm:om1-bc}. Because $\mathcal{B}_1 = 0$ and the constant function lies in $\ker({\mathcal A}_0) \cap \ker(\mathcal{D}_0)$, conditions \eqref{eq:thm2-bc2} and \eqref{eq:thm2-bc3} in Theorem \ref{thm:om1-bc} show that
%\omega_{st}(0) = -\omega_{\varepsilon,0}$. If we multiply the differential equation (\ref{eq:25d}) by $\rho(x)$ and integrate from $0$ to $L$, i.e, take the inner product of it with $1$, and use that $ 1 \in \ker{\mathbf D}({\mathcal A}_0$ and that $\omega_{st}$ weakly satisfies (\ref{eq:25d}), then we obtain
%\begin{equation}
%\label{eq:wave_int}
%  \int_0^L \rho(x) \ddot{\omega}_{st}(x,t) \, \mathrm{d}x = 0 
%\end{equation}
%and thus
%\[
%   \int_0^L \rho(x)\omega_{st}(x,t) \, \mathrm{d}x = \alpha t + \beta.
%\]
%Since $\omega_{st}$ is stable, we must have that $\alpha$ and $\beta$ should be zero. Thus we find that 
\begin{align*}
\omega_1 \int_0^L \rho(x) \mathrm{d}x &= \langle 1, \omega_1 \rangle_X = \mathcal{C}_D(1)  u_0 = \frac{d(L)}{T(L)} u_0, \\
\int_0^L \rho(x) \omega_{0}(x) \mathrm{d}x &= \langle 1, \omega_0 \rangle_{X} = 0. 
\end{align*}
In conclusion, the time domain approach yields
\begin{align}
\begin{split}
  \omega_2(x) =&\ \frac{u_0}{m}, \qquad 
  \omega_1(x) = \frac{d(L)u_0}{m T(L)}, \\
  \omega_0(x) =&\ \left( F(x) - \frac{1}{m} \int_0^L \rho(x) F(x) \, \mathrm{d}x \right) \frac{u_0}{m}. 
  \end{split}
\end{align}
where
\begin{equation}
m = \int_0^L \rho(x) \ \mathrm{d}x, \qquad F(x) = \int_0^x \frac{\int_0^\xi \rho(\eta) \ \mathrm{d}\eta}{T(\xi)} \, \mathrm{d}\xi. \nonumber
\end{equation}
When the material properties $\rho(x) = \rho$, $d(x)=d_0$, and $T(x) = T$ are constant, these formulas reduce to
\begin{equation}
   \omega_2(x) = \frac{u_0}{\rho L}, \ \omega_1(x) = \frac{d_0 u_0}{\rho LT}, \ \omega_0(x) = \frac{3 x^2 - L^2}{6LT} u_0. \label{eq:wave_omegas}
\end{equation}

When the coefficients $\rho(x) = \rho$, $T(x) = T$, and $d(x) = d_0$ are constant, the compliance can also be determined using a frequency domain approach.  Following \cite[Section 10.2]{Ruth} or \cite{curtain2009}, the transfer function $G(x,s)$ from $u(t)$ to $\omega(x,t)$ can be determined from the ansatz $\omega(x,t) = G(x,s)u_0e^{st}$ and $u(t) = u_0 e^{st}$, which gives the equations
\begin{align*}
     &s^2 \rho(x) G(x,s)= 
     (d_0s + T) \frac{\partial^2 G(x,s)}{\partial x^2}, \\
     &\frac{\partial G(0,s)}{\partial x}=0, \quad
     \frac{\partial G(L,s)}{\partial x}=1.
\end{align*}
This leads to the following irrational transfer function
\begin{equation}
    G(x,s)=\frac{\cosh(\tilde{s}x)}{T \tilde{s} \sinh(\tilde{s}L)}, \quad \mbox{ with } \quad  \tilde{s} = s\sqrt{\frac{\rho}{T+d_0s}}.
\end{equation}
Observe that
\begin{equation}
G(x,s) = \frac{1}{TL} \frac{1}{\tilde{s}^2} \cosh(\tilde{s}x) \frac{\tilde{s}L}{\sinh(\tilde{s}L)}. 
\end{equation}
Because for $\xi$ small there holds $\cosh(\xi) = 1 + \xi^2/2 + O(\xi^4)$, $\xi/\sinh(\xi) = 1 - \xi^2/6 + O(\xi^4)$, and $\tilde{s}^2 = s^2 \rho / (T + d_0 s)$, it follows that for $s \approx 0$
\begin{equation}
   G(x,s) = \frac{T+d_0s}{\rho TLs^2} \left(1 + s^2\frac{\rho}{T + d_0 s} \left[\frac{x^2}{2} - \frac{L^2}{6}\right] \right) + O(s^4). \nonumber
\end{equation}
The mass, velocity, and compliance part can now be computed directly from \eqref{eq:G2}--\eqref{eq:G0} as
\begin{align*}
    G_{2}(x) = %\lim_{s\rightarrow 0}s^2 G(x,s) = 
    \frac{1}{\rho L}, \
    G_{1}(x) = %\lim_{s\rightarrow 0}s\left[ G(x,s)-\frac{G_{2}(x)}{s^2}\right] = 
    \frac{d_0}{\rho T L}, \
    G_{0}(x) = %\lim_{s\rightarrow 0}\left[G(x,s)-\frac{G_{2}(x)}{s^2}-\frac{G_1(x)}{s}\right] = 
    \frac{1}{6TL}\left(3x^2 - L^2 \right).
\end{align*}
Looking back at \eqref{eq:wave_omegas}, it is indeed true that $\omega_2(x) = G_2(x) u_0$, 
$\omega_1(x) = G_1(x)u_0$, 
and $\omega_0(x) = G_0(x) u_0$. Observe that the time-domain approach enabled us to compute the compliance with position-dependent coefficients. The frequency domain approach is non-trivial in this situation because there is no closed-form expression for the transfer function.

\subsection{Euler-Bernoulli beam}
\label{Ex-beam}
Now consider the Euler-Bernoulli beam with a damper at one of its endpoints. An Euler-Bernoulli beam with (distributed) structural damping has been studied in \cite{Nikos}. For clarity, it is assumed that the coefficients are constant, although the time domain approach could easily be extended to position-dependent coefficients as well. 

The transversal displacement $\omega(x,t)$ satisfies for $x \in [0,L]$
\begin{align} 
\label{eqp-01}
\rho A_{\mathrm{cs}}&\frac{\partial^2\omega}{\partial t^2}(x,t)%+d \frac{\partial \omega}{\partial t}(x,t)
+EI\frac{\partial^4 \omega}{\partial x^4}(x,t)=0, \\
& EI\frac{\partial^3\omega}{\partial x^3}(0,t)= EI\frac{\partial^2 \omega}{\partial x^2}(0,t) = EI \frac{\partial^2 \omega}{\partial x^2}(L,t) = 0, 
\nonumber \\
-&EI \frac{\partial^3 \omega}{\partial x^3}(L,t) + d\frac{\partial \omega}{\partial t}(L,t) = u(t) , \nonumber
\end{align} 
where the constants $\rho$, $I$, $A_{\mathrm{cs}}$, $E$, $d$, and $u(t)$ are the mass density, second moment of inertia, cross-sectional area, the Young's modulus, and the damping coefficient of damper, and the applied force in the (positive) transversal direction at $x = L$, respectively. This is again a boundary control system of the form \eqref{eq:bc_sys} with state space $X = L^2(0,L)$, inner product
\begin{equation}
  \langle f, g \rangle_X = \rho A_{\mathrm{cs}} \int_0^L f(x) g(x) \ \mathrm{d}x,
\end{equation}
 and
\begin{align*}
\mathcal{A} = \frac{EI}{\rho A_{\mathrm{cs}}} \frac{d^4}{dx^4},
\quad 
%\mathcal{D} = 0, \quad
\mathcal{B}_0 \omega = -EI\frac{d^3 \omega}{dx^3}(L), \quad
\mathcal{B}_1\omega = d\omega(L), \\
\mathbf{D}(\mathcal{A}) = %\mathbf{D}(\mathcal{B}_0) = \mathbf{D}(\mathcal{B}_1) = 
\{ \omega \in H^4(0,L) \mid \tfrac{\partial^3 \omega}{\partial x^3}(0) = \tfrac{\partial^2 \omega}{\partial x^2}(0) = \tfrac{\partial^2 \omega}{\partial x^2}(L) = 0 \},
\end{align*}
$\mathcal{D} = 0$, and $\mathbf{D}(\mathcal{B}_0) = \mathbf{D}(\mathcal{B}_1) = \mathbf{D}(\mathcal{A})$. 
Note that for $f \in \mathbf{D}(\mathcal{A}_0)$ and $g \in \mathbf{D}(\mathcal{A})$
\begin{align*}
\langle f, \mathcal{A}g \rangle &=  EI f(L) \frac{d^3g}{dx^3}(L) + EI \int_0^L \frac{d^2f}{dx^2}(x) \frac{d^2g}{dx^2}(x) \ \mathrm{d}x \\
&= -\mathcal{C}_A(f) \mathcal{B}_0g + \langle \mathcal{A}_0f, g \rangle_X,
\end{align*}
so that $\mathcal{C}_A(f) = f(L)$. As in the previous example this ${\mathcal A}_0$ is also an self-adjoint operaotr.

By noting that the function $20Lx^4 + 12x^5$ is an element of $\mathbf{D}(\mathcal{A})$, it easy to see that $\mathcal{B}_0 : \mathbf{D}(\mathcal{A}) \rightarrow U = \mathbb{R}$ is surjective. By Lemma \ref{La:2}, $\omega_1, \omega_0 \in \mathbf{D}(\mathcal{A})$ and $\omega_0$, $\omega_1$, and $\omega_2$ are strong solutions of \eqref{eq:25a}, \eqref{eq:25b}, and \eqref{eq:25c}.  

Note that $\mathcal{A} \omega_2 = 0$ implies that $\omega_2(x)$ is of the form
\begin{equation}
\omega_2(x) = C_1 + C_2x + C_3 x^2 + C_4 x^3,
\end{equation}
for some constants $C_1$, $C_2$, $C_3$, and $C_4$. Because $\omega_2 \in \mathbf{D}(\mathcal{A}_0)$, $C_3 = C_4 = 0$. 

Because $\mathcal{D} = 0$, \eqref{eq:25b} reduces to $\mathcal{A}\omega_1 = 0$ and $\mathcal{B}_0\omega_1 + \mathcal{B}_1 \omega_2 = 0$.  The requirement that $\mathcal{A}\omega_1 = 0$ and $\omega_1 \in \mathbf{D}(\mathcal{A})$ implies similarly as above that
\begin{equation}
\omega_1(x) = C_5 + C_6 x,  \label{eq:beam_om1}
\end{equation}
for certain coefficients $C_5$ and $C_6$. 
Note that $\mathcal{B}_0 \omega_1 = 0$ for all $C_5$ and $C_6$, so that the requirement $\mathcal{B}_0 \omega_1 + \mathcal{B}_1 \omega_2 = 0$ reduces to $\mathcal{B}_1 \omega_2 = 0$ which implies that $C_1 + C_2L = 0$.   

To determine $\omega_0$, note that because $\mathcal{D} = 0$, \eqref{eq:25c} reduces to $\omega_2 + \mathcal{A}\omega_0 = 0$ and $\mathcal{B}_0 \omega_0 + \mathcal{B}_1\omega_1 = u_0$. The equation ${\mathcal A}\omega_0 = - \omega_2$  together with the zero boundary conditions on $\omega_0$ at $x = 0$, gives that %\marginpar{I got a minus}
\begin{equation}
\omega_0(x) = -\frac{\rho A_{cs}}{EI}\left(\frac{C_1}{24} x^4 + \frac{C_2}{120 } x^5  \right) + C_{7} + C_{8}x, \label{eq:beam_om0}
\end{equation}
for some constants $C_7$ and $C_8$. The (boundary) condition $\mathcal{B}_0 \omega_0 + \mathcal{B}_1\omega_1 = u_0$ implies that
\begin{equation}
  \rho A_{\mathrm{cs}}(C_1L + \tfrac{1}{2}C_2 L^2) + d(C_5 + C_6L) = u_0. \label{eq:beam_bc0}
\end{equation}

Because $L-x$ is an element of $\ker(\mathcal{A}_0) \cap \ker(\mathcal{D}_0) \cap \ker(\mathcal{C}_A)$, \eqref{eq:thm2-bc1-a} in Theorem \ref{thm:om1-bc} shows that
\begin{multline*}
\rho A_{\mathrm{cs}} L^2\left(\tfrac{1}{2}C_1 + \tfrac{1}{6}C_2L \right) = \langle L-x, \omega_2 \rangle_X \\
= \mathcal{C}_A (L-x) (u_0 - \mathcal{B}_1 \omega_1)  = 0, 
\end{multline*}
which together with $C_1 + C_2L = 0$ implies that $C_1 = C_2 = 0$, and thus that $\omega_2(x) \equiv 0$. 

Calculating \eqref{eq:thm2-bc2-a} and \eqref{eq:thm2-bc3-a} with $\phi = L-x$ shows that
\begin{align}
\rho A_{\mathrm{cs}}L^2 \left( \tfrac{1}{2}C_5 + \tfrac{1}{6}C_6L \right)=\langle L-x, \omega_1 \rangle_X = 0, 
\label{eq:beam_C56} \\
\rho A_{\mathrm{cs}}L^2 \left( \tfrac{1}{2}C_7 + \tfrac{1}{6}C_8L \right)=\langle L-x, \omega_0 \rangle_X = 0. \label{eq:beam_C78}
\end{align}
Equation \eqref{eq:beam_C56} together with \eqref{eq:beam_bc0} shows that $C_6 = 3 u_0 / (2Ld)$ and $C_5 = -u_0/(2d)$. 

Finally, note that the constant function is an element of $\ker(\mathcal{A}_0) \cap \ker(\mathcal{D}_0)$ and that \eqref{eq:thm2-bc2-a} thus shows that
\begin{equation}
\frac{\rho A_{\mathrm{cs}}Lu_0}{4d} = \langle 1, \omega_1 \rangle_X = \mathcal{C}_A(1) \mathcal{B}_1 \omega_0 = d (C_7 + C_8 L)
\end{equation}
which, together with \eqref{eq:beam_C78}, shows that $C_7 = \rho A_{\mathrm{cs}}L / (8d^2)$ and $C_8 = -3\rho A_{\mathrm{cs}} / (8d^2)$. The time-domain approach thus yields
\begin{equation}
\omega_2(x) = 0, \
\omega_1(x) = \frac{3x-L}{2Ld}, \
\omega_0(x) = \rho A_{\mathrm{cs}}\frac{L - 3x}{8d^2}. \label{eq:beam_oms}
\end{equation}

These results are again validated by a frequency domain approach, which is rather cumbersome compared to the time domain approach considered before. Following \cite[Section 10.2]{Ruth} or \cite{curtain2009}, the transfer function can be determined from the ansatz $\omega(x,t) = G(x,s)u_0e^{st}$ and $u(t) = u_0 e^{st}$, which gives
\begin{align}
     \rho A_{\mathrm{cs}} s^2 & G(x,s)+ 
     EI \frac{\partial^4 G(x,s)}{\partial x^4} = 0, \label{eq:Gbeam_pde} \\
     &\frac{\partial^2 G(0,s)}{\partial x^2}=\frac{\partial^3 G(0,s)}{\partial x^3}=\frac{\partial^2 G(L,s)}{\partial x^2}=0, \nonumber \\
     &-EI\frac{\partial^3 G(L,s)}{\partial x^3} + d s G(L,s)=1. \label{eq:Gbeam_bc}
\end{align}
The solution of \eqref{eq:Gbeam_pde} is of the form
\begin{equation}
G(x,s) = 
\begin{bmatrix}
C_1(s) & C_2(s) & C_3(s) & C_4(s)
\end{bmatrix} \mathbf{b}(x,s), \label{eq:Gbeam_form}
\end{equation}
for certain constants $C_1(s)$, $C_2(s)$, $C_3(s)$, and $C_4(s)$ that do not depend on $x$ and with
\begin{equation}
\mathbf{b}(x,s) = 
%\begin{bmatrix}
%b_1(x,s) \\ b_2(x,s) \\ b_3(x,s) \\ b_4(x,s)
%\end{bmatrix} = 
\begin{bmatrix}
\cosh(\beta(s)x)\cos(\beta(s)x) \\
\sinh(\beta(s)x)\cos(\beta(s)x) \\
\cosh(\beta(s)x)\sin(\beta(s)x) \\\sinh(\beta(s)x)\sin(\beta(s)x)
\end{bmatrix}, \
\beta(s) = \sqrt[4]{\frac{\rho A_{\mathrm{cs}}s^2}{4EI}}. \nonumber
\end{equation} 
The coefficients $C_i(s)$ are determined by substituting \eqref{eq:Gbeam_form} into the boundary conditions \eqref{eq:Gbeam_bc}. 
Because computations now becomes extremely cumbersome by hand, the Taylor series of $s^2G(x,s)$ around $s = 0$ is computed using the symbolic toolbox in GNU Octave. This yields
\begin{equation}
s^2 G(x,s) = s\frac{3x-L}{2Ld} + s^2 \rho A_{\mathrm{cs}}  \frac{(L-3x)}{8d^2} +  O(s^3).
\end{equation}
It now easy to verify from this expression and the expressions for $\omega_2(x)$, $\omega_1(x)$, and $\omega_0(x)$ in \eqref{eq:beam_oms} that $\omega_2(x) = G_2(x) u_0 = 0$, $\omega_1(x) = G_1(x) u_0$, and $\omega_0(x) = G_0(x) u_0$. Also observe that the time-domain approach was again much easier to apply in this example than the frequency domain approach because expanding $s^2G(s,x)$ in a Taylor series around $s = 0$ by hand is very cumbersome.

\subsection{Kirchhoff plate model}
\label{example}
The proposed approach is also applied to the rectangular Kirchhoff plate model. The main assumption is that the plate structure fibers which are orthogonal to the middle plane remain orthogonal after deformation. The deformation of the plate is represented by the transversal displacement of the midplane $\omega(x,y,t)$ which is defined on $\Omega = (0,L) \times (0,W)$. The input $u(t)$ is a uniform force applied to the half edge $(x,y) \in \{ x \} \times [0, W/2]$ and the other edges are stress free. 

This model can be written in the general form of a boundary control system, in which the operator ${\mathcal A}$ on the  state space $X = L^2(\Omega)$ is given by
\begin{equation}
  {\mathcal A}\omega = \tilde{D}\left( \frac{\partial^4 \omega}{\partial x^4}(x,y,t) + 2 \frac{\partial^4 \omega}{\partial x^2 \partial y^2}(x,y,t) + \frac{\partial^4 \omega}{\partial y^4}(x,y,t)\right), \nonumber
\end{equation}
with $\tilde{D} = Eh^2/(12\rho(1-\nu^2))$ where $E$ is the Young's modulus, $\nu \in [0,\tfrac{1}{2}]$ is Poisson's ratio, $h$ is the thickness of the plate, and $\rho$ is the mass density. The domain $\mathbf{D}(\mathcal{A})$ are the functions in $H^4(\Omega)$ satisfying the following stress free boundary conditions, see, e.g., \cite[Chapter 13]{Zienkiewicz},
\begin{align}
\label{plate_BC1}
&M_{x} = M_{xy} = 0, \qquad\qquad\qquad\quad \mathrm{on}\ \{ 0, a \} \times (0,W), \\
\label{plate_BC2}
&Q_{x}=0, \qquad\qquad \mathrm{on}\ \{0 \} \times (0,W/2) \cup \{ a \} \times (0,W),\\
\label{plate_BC3}
&M_{y} = M_{xy} = Q_y = 0, \qquad\qquad\quad \mathrm{on}\ (0,L) \times \{0,b \},
\end{align}
where
\begin{align*}
M_x = -\tilde{D} \left( \frac{\partial^2 \omega}{\partial x^2} + \nu \frac{\partial^2 \omega}{\partial y^2}\right), \
M_y = -\tilde{D}\left( \frac{\partial^2 \omega}{\partial y^2} + \nu \frac{\partial^2 \omega}{\partial x^2}\right),
\\
M_{xy} = -\tilde{D}(1-\nu) \frac{\partial^2 \omega}{\partial x \partial y}, \qquad\qquad\qquad
\\
Q_x = \tilde{D}\left( \frac{\partial^3 \omega}{\partial x^3} + \frac{\partial^3 \omega}{\partial x \partial y^2}\right), \
Q_y = \tilde{D}\left( \frac{\partial^3 \omega}{\partial y^3} + \frac{\partial^3 \omega}{\partial y \partial x^2}\right).
\end{align*}
The plate is structurally damped, so for some $\beta > 0$
\begin{equation}
\mathcal{D} = \beta \mathcal{A}, \qquad \qquad \mathbf{D}(\mathcal{D}) = \mathbf{D}(\mathcal{A}).
\end{equation}

The control space is chosen as $U = L^2(W/2, W)$ and the domain $\mathbf{D}(\mathcal{B}_0)$ consists of all functions in $H^{3+1/2}(\Omega)$ that satisfy the boundary conditions in \eqref{plate_BC1}--\eqref{plate_BC3}, and for $y \in (W/2, W)$
\begin{equation}
\mathcal{B}_0 \omega = \tilde{D}\left( \frac{\partial^3 \omega}{\partial x^3}(0,y) + \frac{\partial^3 \omega}{\partial x \partial y^2}(0,y)\right). 
\end{equation}
By the trace theorem, $\mathcal{B}_0(\mathbf{D}(\mathcal{A})) \subset \mathcal{B}_0(H^4(\Omega)) = H^{1/2}(W/2,W)$. It thus follows that $\mathcal{B}_0 : \mathbf{D}(\mathcal{A}) \rightarrow U$ is not surjective and that Lemma \ref{La:2} therefore does not apply. In particular, it is now only known that $\omega_1$ and $\omega_0$ are weak solutions of \eqref{eq:25b-w} and \eqref{eq:25c-w}, but they are generally not strong solutions of \eqref{eq:25b} and \eqref{eq:25c}. In the following only inputs $u \in U$ that are constant over $[W/2, W]$ will be considered. However, because $\partial^3 \omega / \partial x^3 + \partial^3 \omega / \partial x \partial y^2$ contains a jump at $y = W/2$, these inputs do not lie in the range of $\mathcal{B}_0 : \mathbf{D}(\mathcal{A}) \rightarrow U$.

Note that \eqref{eq:25a} shows that $\omega_2$ is in the kernel $\mathcal{A}_0$ and that a standard computation using integration by parts shows that for every $\omega \in \mathbf{D}(\mathcal{A}_0)$
\begin{multline}
\langle \omega, \mathcal{A}_0 \omega \rangle = \\ \tilde{D}\iint_\Omega \begin{bmatrix}
\frac{\partial^2 \omega}{\partial x^2} \\
\frac{\partial^2 \omega}{\partial y^2} \\ 
\frac{\partial^2 \omega}{\partial x \partial y}
\end{bmatrix}^\top \begin{bmatrix}
1 & \nu & 0 \\
\nu & 1 & 0 \\
0 & 0 & 2(1-\nu)
\end{bmatrix}
\begin{bmatrix}
\frac{\partial^2 \omega}{\partial x^2} \\
\frac{\partial^2 \omega}{\partial y^2} \\ 
\frac{\partial^2 \omega}{\partial x \partial y}
\end{bmatrix} \ \mathrm{d}x \ \mathrm{d}y. \label{eq:plate_weak3}
\end{multline}
Because the matrix in \eqref{eq:plate_weak3} is positive definite for $\nu \in [0,\tfrac{1}{2}]$, it follows that $\mathcal{A}_0\omega_2 = 0$ implies that
\begin{equation}
\frac{\partial^2 \omega_2}{\partial x^2}  = 0, \qquad\quad 
\frac{\partial^2 \omega_2}{\partial y^2}  = 0, \qquad\quad
\frac{\partial^2 \omega_2}{\partial x \partial y} = 0. \label{eq:plate_kern_cond}
\end{equation}
It is now easy to see that this implies that $\omega_2$ is of the form
\begin{equation}
\omega_2(x,y) = a_1 + a_2 x + a_3 y, \label{eq:plate_kernel}
\end{equation}
for some constants $a_1, a_2, a_3 \in \mathbb{R}$. 
Note that the basis function $1$ represents the translation of the plate in the out-of-plane direction, $x$ represents the rotation around the $x = 0$ axis, and $y$ represents the rotation around the $y = 0$ axis. 

The remainder of the computations cannot be carried out analytically. This is also related to the observation from before that $\omega_1$ and $\omega_0$ are not strong solutions of \eqref{eq:25b} and \eqref{eq:25c} because $\mathcal{B}_0 : \mathbf{D}(\mathcal{A}) \rightarrow U$ is not surjective. The plate equation is therefore discretized with the finite element scheme proposed by Bogner et al.\ \cite{Bogner}, see also \cite[Chapter 13]{Zienkiewicz}. In this scheme, each rectangular element was 4 nodes at its four corners. Each node has four degrees of freedom: $\omega$, $\partial \omega / \partial x$, $\partial \omega / \partial y$, and $\partial^2 \omega / \partial x \partial y$. The used parameter values are given in Table \ref{tab:plate}. There are $M_x = 30$ elements in the $x$-direction and $M_y = 40$ in the $y$-direction. This leads to a (row)vector $\mathbf{N}(x,y) \in \mathbb{R}^{1\times N}$ (with $N = 4(M_x+1)(M_y+1) = 5084$) of finite-element shape functions and the approximation
\begin{equation}
\omega(x,y,t) \approx \mathbf{N}(x,y) \mathbf{w}(t),
\end{equation}
where $\mathbf{w}(t) \in \mathbb{R}^N$ contains the nodal values of $\omega$, $\partial \omega / \partial x$, $\partial \omega / \partial y$, and $\partial^2 \omega / \partial x \partial y$. A Galerkin discretization for this set of shape functions yields a system of ODEs, see (\ref{eq:bc_sys-weak}),
\begin{equation}
\mathbf{M} \ddot{\mathbf{w}}(t) + \beta \mathbf{K}\dot{\mathbf{w}}(t) + \mathbf{K} \mathbf{w}(t) = \mathbf{B}u(t) + \beta \mathbf{B} \dot{u}(t), \label{eq:plate_dyn}
\end{equation}
where $\mathbf{M}, \mathbf{K} \in \mathbb{R}^{N \times N}$ are given by
\begin{align*}
\mathbf{M} &= \rho h \iint_\Omega \mathbf{N}^\top \mathbf{N} \ \mathrm{d}x \ \mathrm{d}y, \\ \mathbf{K} &= \tilde{D}\iint_\Omega \begin{bmatrix}
\frac{\partial^2 \mathbf{N}^\top}{\partial x^2} \\
\frac{\partial^2 \mathbf{N}^\top}{\partial y^2} \\
\frac{\partial^2 \mathbf{N}^\top}{\partial x \partial y}
\end{bmatrix} \begin{bmatrix}
1 & \nu & 0 \\
\nu & 1 & 0 \\
0 & 0 & 2(1-\nu)
\end{bmatrix} \begin{bmatrix}
\frac{\partial^2 \mathbf{N}^\top}{\partial x^2} \\
\frac{\partial^2 \mathbf{N}^\top}{\partial y^2} \\
\frac{\partial^2 \mathbf{N}^\top}{\partial x \partial y}
\end{bmatrix}^\top \ \mathrm{d}x \mathrm{d}y. 
\end{align*}
The input matrix $\mathbf{B} \in \mathbb{R}^{N\times 1}$ contains 1's at the locations corresponding to the transversal displacement of the nodes on the edge $(x,y) \in \{ 0\} \times [W/2, W]$ and zeros otherwise. Note that the appearance of $\dot{u}(t)$ in \eqref{eq:plate_dyn} is similar as in the weak form of general boundary control system \eqref{eq:bc_sys-weak}, but that \eqref{eq:plate_dyn} cannot be obtained directly from \eqref{eq:bc_sys-weak} by restricting $\phi$ and $\omega(t)$ to the space spanned by the FE shape functions.

\begin{table}
\caption{Parameter values used in the plate example}
\label{tab:plate}
\centering
\begin{tabular}{|l|c|c|l|c|}
\cline{1-2} \cline{4-5}
{\bf Parameter} & {\bf Value} & \multicolumn{1}{|c|}{} & {\bf Parameter} & {\bf Value} \\ \cline{1-2} \cline{4-5}
$\rho$ & 2700 kg/m$^3$ & & $h$ & 0.01 m \\ \cline{1-2} \cline{4-5}
$E$    & 69 GPa & & $L$ & 0.5 m \\ \cline{1-2} \cline{4-5}
$\nu$  & 0.3 & & $W$    & 0.4 m \\ \cline{1-2} \cline{4-5}
$\beta$ & 0.01 s \\ \cline{1-2}
\end{tabular}

\end{table}

\begin{figure}[htb]
\centering
\begin{subfigure}[b]{0.47\textwidth}
\centering
\includegraphics[width=\textwidth]{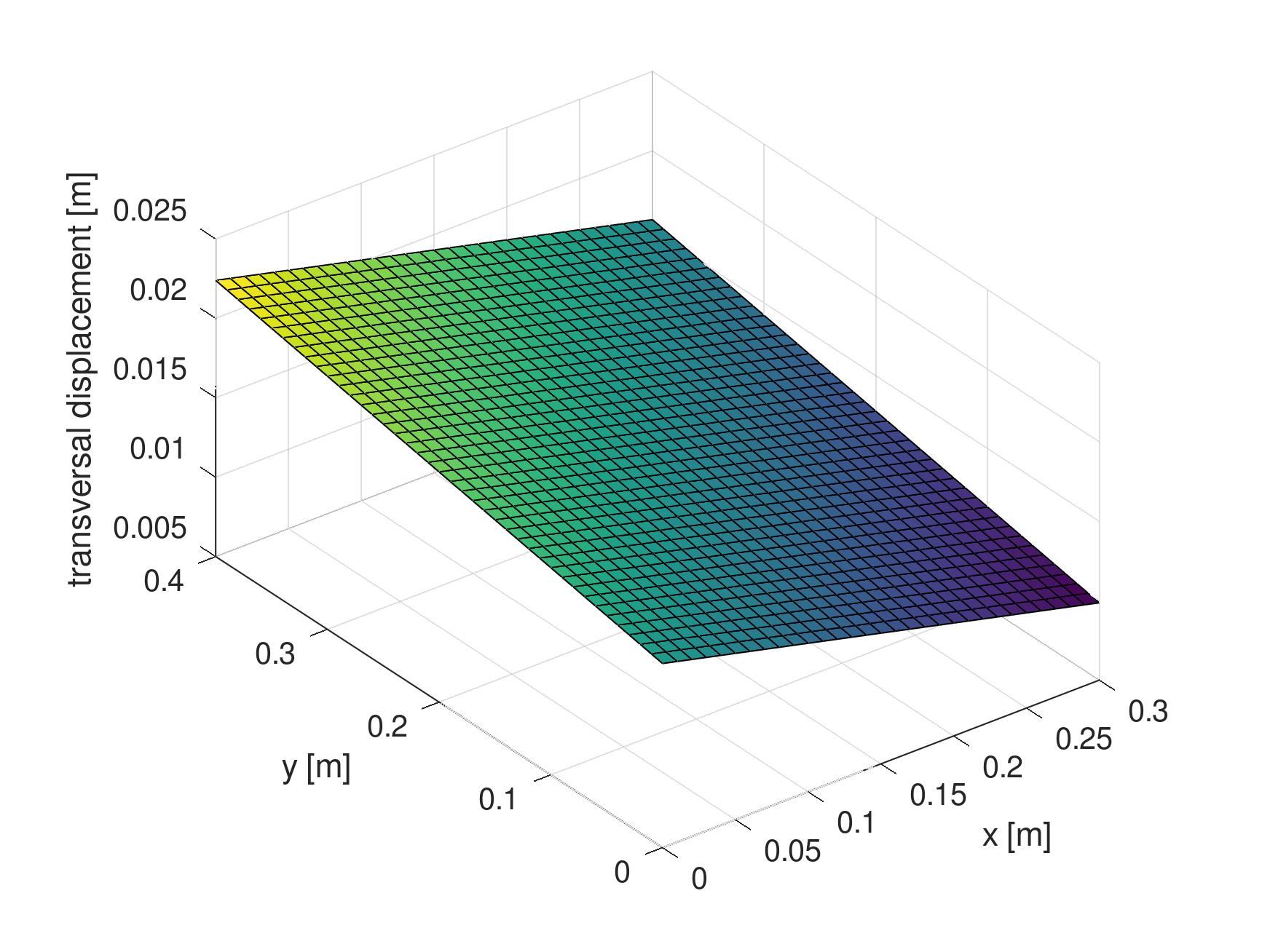}
\caption{$\mathbf{N}(x,y)\mathbf{w}_2 \approx \omega_2(x,y)$}
\label{fig:plate_omega2}
\end{subfigure}
\begin{subfigure}[b]{0.47\textwidth}
\centering
\includegraphics[width=\textwidth]{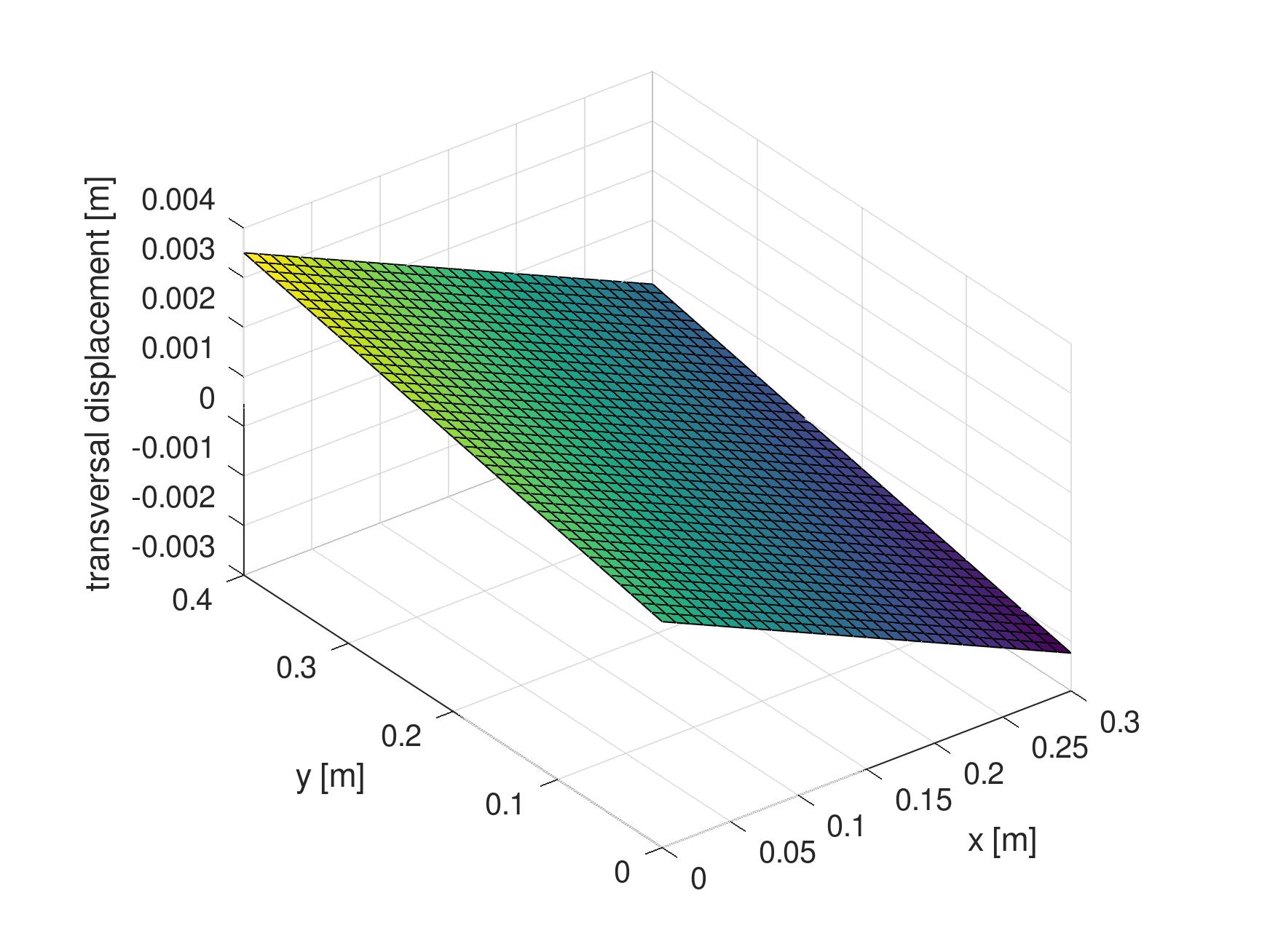}
\caption{$\mathbf{N}(x,y)\mathbf{w}_1 \approx \omega_1(x,y)$}
\label{fig:plate_omega1}
\end{subfigure}
\begin{subfigure}[b]{0.47\textwidth}
\centering
\includegraphics[width=\textwidth]{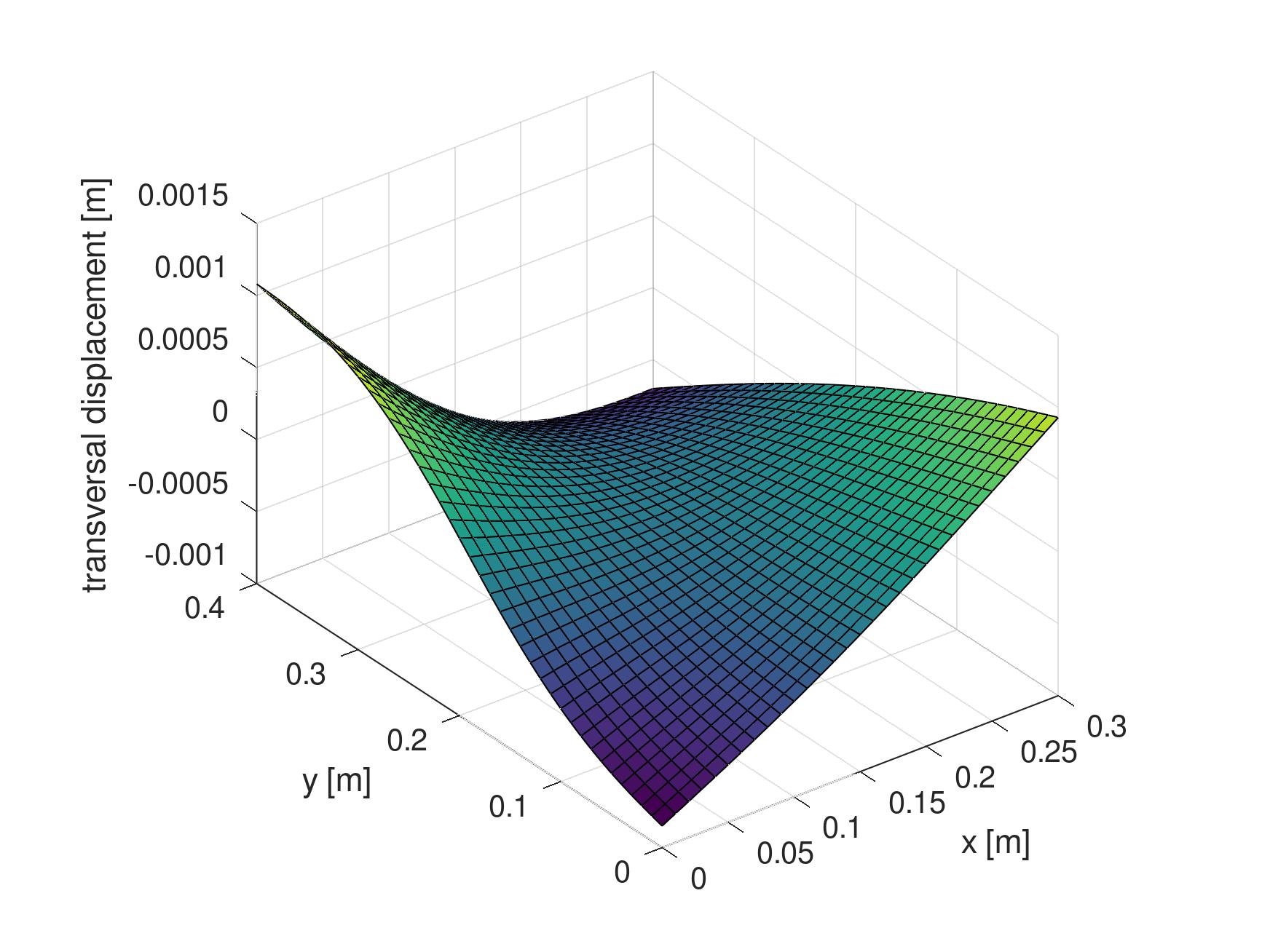}
\caption{$\mathbf{N}(x,y)\mathbf{w}_0 \approx \omega_0(x,y)$}
\label{fig:plate_omega0}
\end{subfigure}
\caption{The obtained $\mathbf{N}(x,y)\mathbf{w}_2$, $\mathbf{N}(x,y)\mathbf{w}_1$, and $\mathbf{N}(x,y) \mathbf{w}_0$ for the plate with structural damping when a force in the positive $z$-direction is applied on the edge $(x,y) \in \{ 0\} \times [W/2, W]$. }
\label{fig:plate}
\end{figure}

Recall that the kernel of ${\mathcal A}_0$ is spanned by $\{1,x,y \}$, see \eqref{eq:plate_kernel}. These functions can be represented exactly on the finite element basis, i.e.\ there exist $\mathbf{v}_1, \mathbf{v}_x, \mathbf{v}_y \in \mathbb{R}^N$ such that
\begin{equation}
1 = \mathbf{N}(x,y)\mathbf{v}_1, \quad 
x = \mathbf{N}(x,y)\mathbf{v}_x, \quad
y = \mathbf{N}(x,y)\mathbf{v}_y.
\end{equation}
It is now easy to see that the kernel of $\mathbf{K}$ is spanned by $\mathbf{v}_1$, $\mathbf{v}_x$, and $\mathbf{v}_y$, i.e.\ by the (column) range of the matrix $\mathbf{V}  \in \mathbb{R}^{N\times 3}$
\begin{equation}
\mathbf{V} = \begin{bmatrix}
\mathbf{v}_1 & \mathbf{v}_x & \mathbf{v}_y
\end{bmatrix}.
\end{equation}

Note that the finite element discretization has transformed the boundary control system into a system with internal control on a finite-dimensional space of the form \eqref{eqm-1} (after multiplying \eqref{eq:plate_dyn} by $\mathbf{M}^{-1}$ from the left). By Theorem \ref{thm:om1}, $\mathbf{w}_1$ and $\mathbf{w}_2$ lie in the kernel of $\mathbf{K}$ and can thus be represented as $\mathbf{w}_2 = \mathbf{V} \mathbf{a}_2$ and $\mathbf{w}_1 = \mathbf{V}\mathbf{a}_1$. Inserting these expressions into \eqref{eq:8c}, it is now easy to see that $\mathbf{a}_2$ and $\mathbf{w}_0$ are the unique solution of the linear system
\begin{equation}
\begin{bmatrix}
\mathbf{K} & \mathbf{M}\mathbf{V} \\
\mathbf{V}^\top\mathbf{M} & 0
\end{bmatrix} \begin{bmatrix}
\mathbf{w}_{0} \\ \mathbf{a}_2
\end{bmatrix} =  \begin{bmatrix}
\mathbf{B} \\ 0
\end{bmatrix}, \label{eq:w0a2}
\end{equation}
in which the second line represents the first condition in \eqref{eq:thm2}. The second condition in \eqref{eq:thm2} (with the inner product $\langle \mathbf{v}, \mathbf{w} \rangle := \mathbf{v}^\top \mathbf{M} \mathbf{w}$) shows that 
\[ \mathbf{a}_1 = (\mathbf{V}^\top \mathbf{M} \mathbf{V})^{-1} \mathbf{V}^\top \mathbf{B}_1 u_0. \]
The $\mathbf{w}_2 = \mathbf{V}\mathbf{a}_2$, $\mathbf{w}_1 = \mathbf{V} \mathbf{a}_1$, and $\mathbf{w}_0$ obtained in this way are displayed in Figure \ref{fig:plate}. 

\begin{figure}[htb]
\centering
\begin{subfigure}[b]{0.45\textwidth}
\centering
\includegraphics[width=\textwidth]{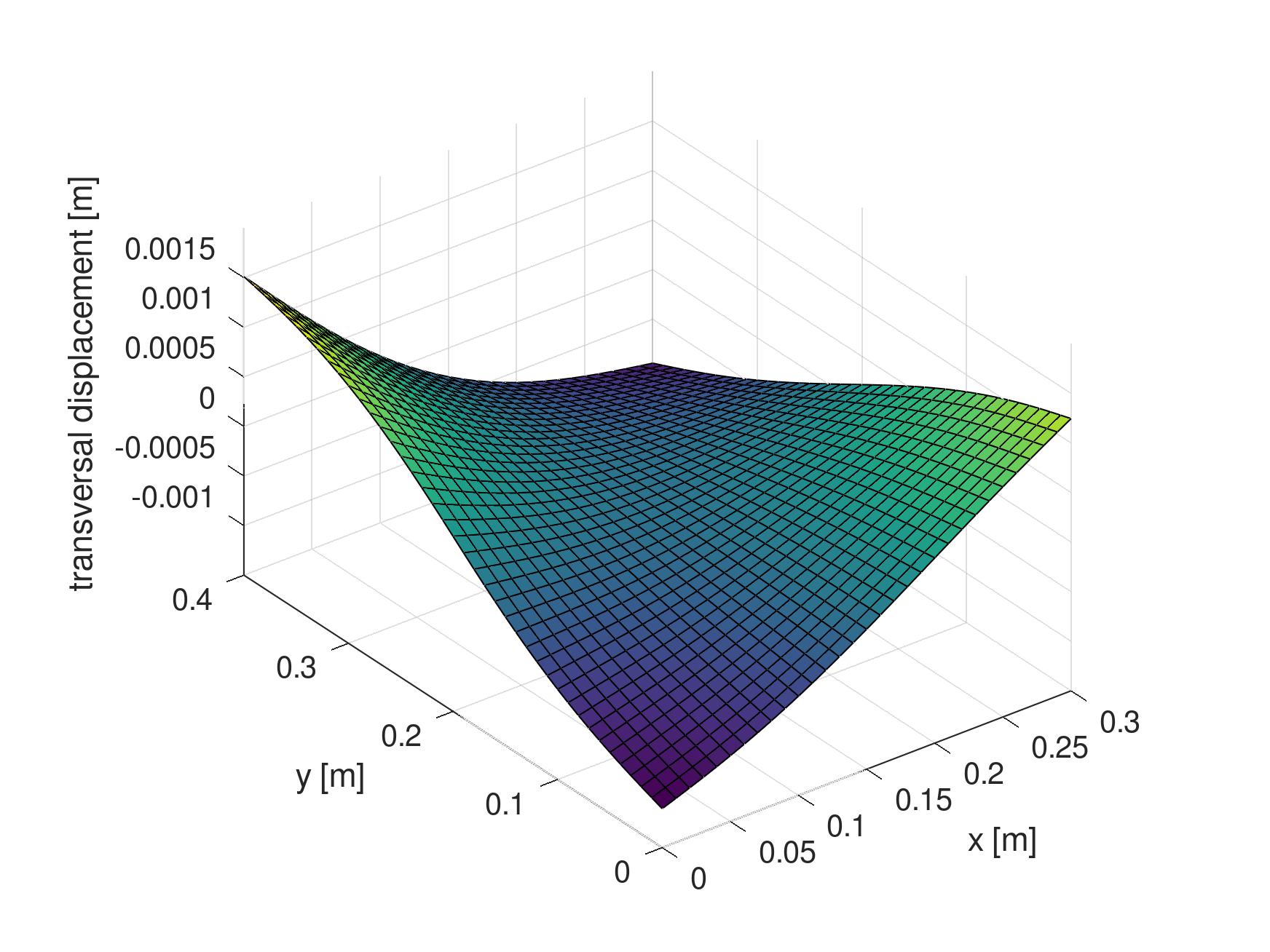}
\caption{$t = 0.25$ s}
\end{subfigure}
\begin{subfigure}[b]{0.45\textwidth}
\centering
\includegraphics[width=\textwidth]{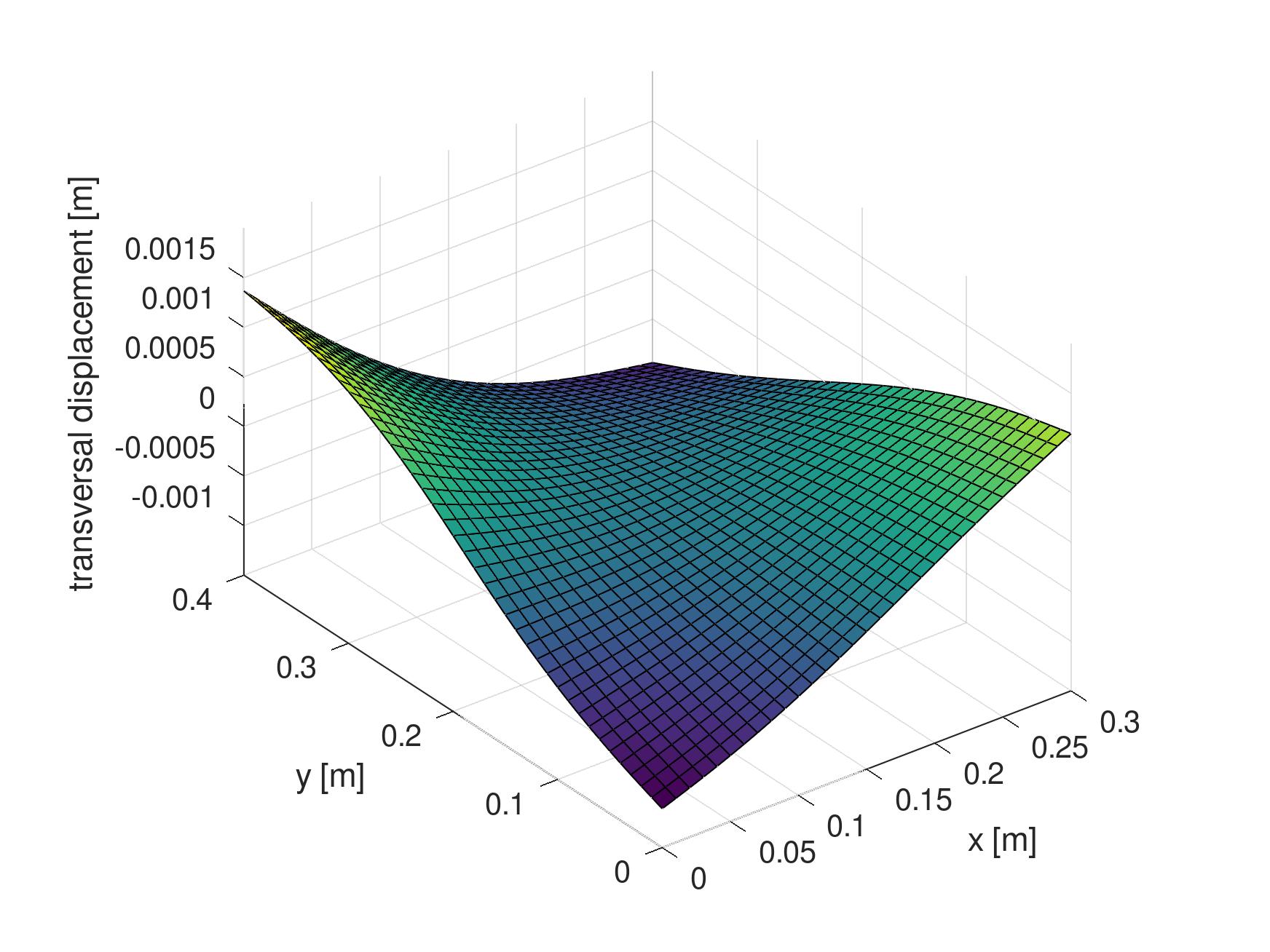}
\caption{$t = 0.5$ s}
\end{subfigure}
\begin{subfigure}[b]{0.45\textwidth}
\centering
\includegraphics[width=\textwidth]{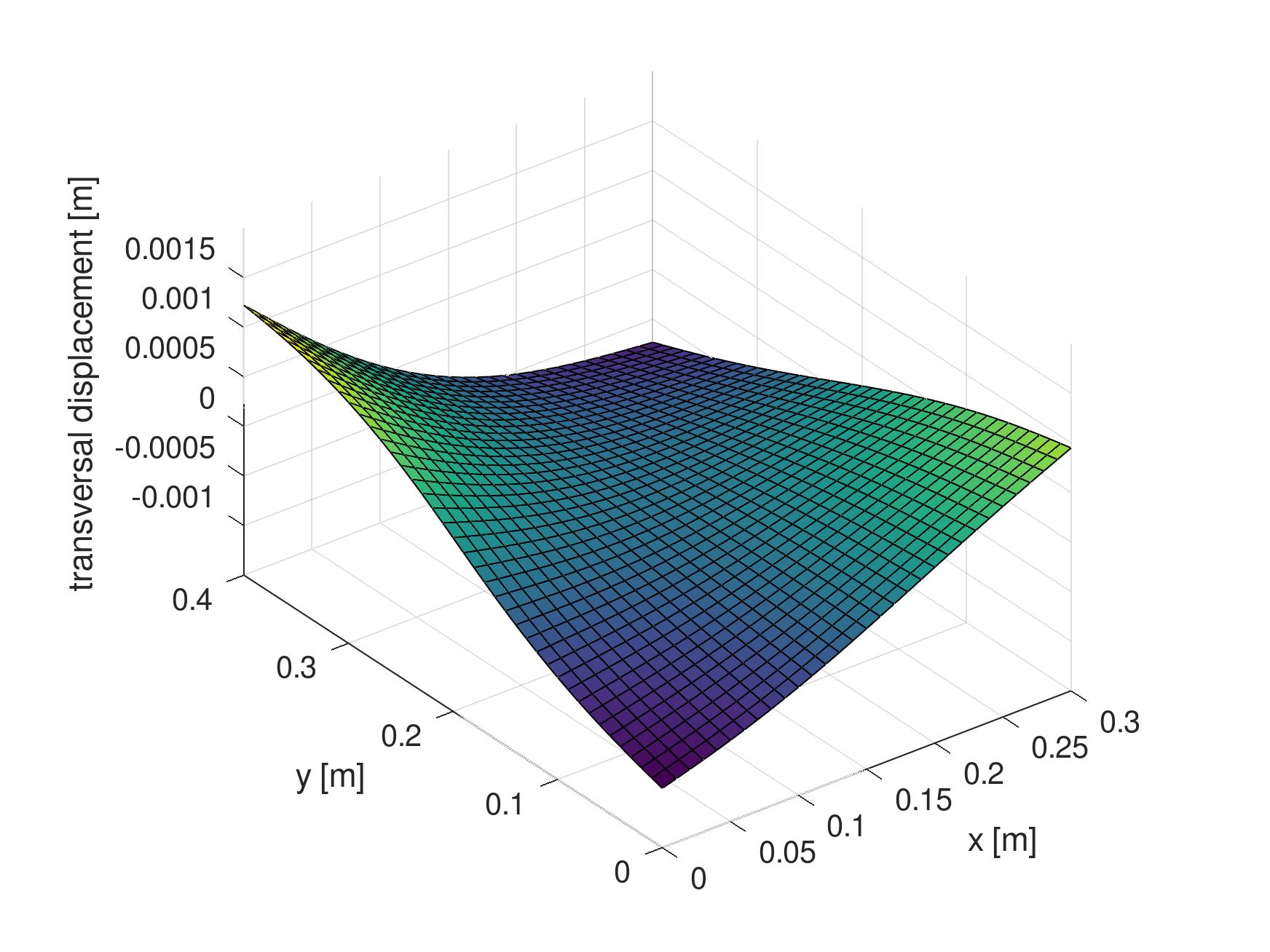}
\caption{$t = 1$ s}
\end{subfigure}
\begin{subfigure}[b]{0.45\textwidth}
\centering
\includegraphics[width=\textwidth]{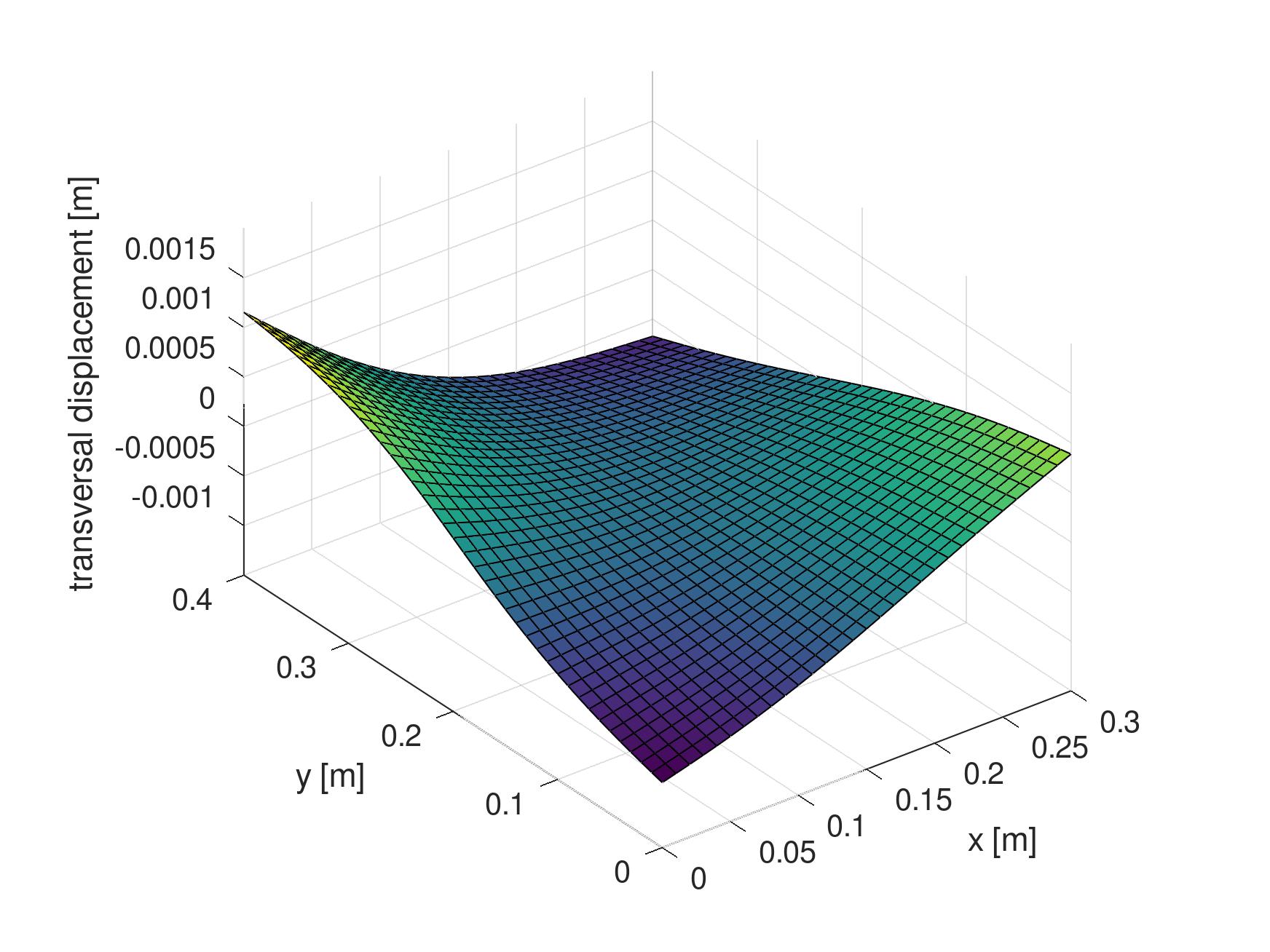}
\caption{$t = 4$ s}
\end{subfigure}
\caption{Four snapshots of $\mathbf{N}(x,y) \mathbf{w}(t) - \tfrac{1}{2}t^2\mathbf{N}(s,t)\mathbf{w}_2 - t\mathbf{w}_1$ of the simulation of the plate with structural damping. A unit force in the positive $z$-direction is applied on the edge $(x,y) \in \{ 0\} \times [W/2, W]$. }
\label{fig:plate_snapshots_omega0}
\end{figure}

To validate the obtained results, \eqref{eq:plate_dyn} with $u(t) \equiv u_0 = 1$ is integrated over time and it is verified that $\mathbf{w}_{st}(t) - \tfrac{1}{2}t^2 \mathbf{w}_2 - t \mathbf{w}_1 \rightarrow \mathbf{w}_0$ for $t \rightarrow \infty$. Time integration is done by writing \eqref{eq:plate_dyn} as a first order ODE
\begin{equation}
\frac{d}{dt}
\begin{bmatrix}
\mathbf{M} & 0 \\
0 & \mathbf{M}
\end{bmatrix} \begin{bmatrix}
\mathbf{w}(t) \\ \dot{\mathbf{w}}(t)
\end{bmatrix} = 
\begin{bmatrix}
0 & \mathbf{M} \\
-\mathbf{K} & -\mathbf{D}
\end{bmatrix} \begin{bmatrix}
\mathbf{w}(t) \\ \dot{\mathbf{w}}(t)
\end{bmatrix},
\end{equation}
with initial conditions $\mathbf{w}(0) = \mathbf{0}$ and $ \dot{\mathbf{w}}(0) = \beta \mathbf{B}u_0$ and applying the Crank-Nicholson scheme with a fixed time step of $\Delta t = 0.001$. %Figure \ref{fig:plate_snapshots} shows four snapshots of $\mathbf{N}(x,y)\mathbf{w}(t)$. The influence of the compliance $\mathbf{w}_0$ is negligible compared the rigid body motion proportional to $\mathbf{N}(x,y) \mathbf{w}_2$ in Figure \ref{fig:plate_omega2}.
Figure \ref{fig:plate_snapshots_omega0} shows four snapshots of $\mathbf{N}(x,y)\mathbf{w}(t) - \tfrac{1}{2}t^2 \mathbf{N}(x,y) \mathbf{w}_2 - t \mathbf{w}_1$. The obtained shape is very close to the compliance $\mathbf{N}(x,y) \mathbf{w}_0$ in Figure \ref{fig:plate_omega0}. Indeed, the maximal difference between the transversal displacement components of $\mathbf{w}_0$ and $\mathbf{w}(4)$ is less then $2.3 \cdot 10^{-6}$ (0.6 \%). 

It is also possible to approximate $\mathbf{w}_0$, $\mathbf{w}_1$, and $\mathbf{w}_2$ based on the simulation of $\mathbf{w}(t)$ only using the approach described at the end of Section \ref{sec:2}. This leads to approximations $\tilde{\mathbf{w}}_0$, $\tilde{\mathbf{w}}_1$, and $\tilde{\mathbf{w}}_2$ as in \eqref{eq:omegai_sol}. Between $t = 3$ and $t = 4$, $n = 11$ time instances with a uniform step size of $0.1$ are considered. The difference between the displacement components of the obtained $\tilde{\mathbf{w}}_0$, $\tilde{\mathbf{w}}_1$, and $\tilde{\mathbf{w}}_2$ and the $\mathbf{w}_0$, $\mathbf{w}_1$, and $\mathbf{w}_2$ displayed in Figure \ref{fig:plate} is below $5 \cdot 10^{-5}$. The relative errors in the displacement components of $\mathbf{w}_2$, $\mathbf{w}_1$, and $\mathbf{w}_0$ are $0.003 \%$, $0.9 \%$, and $5 \%$, respectively. The obtained $\tilde{\mathbf{w}}_0$, $\tilde{\mathbf{w}}_1$, and $\tilde{\mathbf{w}}_2$ are thus visually indistinguishable from the $\mathbf{w}_0$, $\mathbf{w}_1$, and $\mathbf{w}_2$ displayed in Figure \ref{fig:plate}.

%\begin{figure}
%\centering
%\begin{subfigure}[b]{0.35\textwidth}
%\centering
%\includegraphics[width=\textwidth]{Fig_plate_snapshot1}
%\caption{$t = 0.25$ s}
%\end{subfigure}
%\begin{subfigure}[b]{0.35\textwidth}
%\centering
%\includegraphics[width=\textwidth]{Fig_plate_snapshot2}
%\caption{$t = 0.5$ s}
%\end{subfigure}
%\begin{subfigure}[b]{0.35\textwidth}
%\centering
%\includegraphics[width=\textwidth]{Fig_plate_snapshot3}
%\caption{$t = 1$ s}
%\end{subfigure}
%\begin{subfigure}[b]{0.35\textwidth}
%\centering
%\includegraphics[width=\textwidth]{Fig_plate_snapshot4}
%\caption{$t = 4$ s}
%\end{subfigure}
%\caption{Four snapshots of the simulation of the plate with structural damping. A unit force in the positive $z$-direction is applied on the edge $(x,y) \in \{ 0\} \times [W/2, W]$. }
%\label{fig:plate_snapshots}
%\end{figure}

\begin{thebibliography}{99}

\bibitem{Ruth}
Curtain, R., Zwart, H., {\em Introduction to Infinite-Dimensional Systems Theory:
A State-Space Approach}, Springer Verlag, New-York, 2020.

\bibitem{Fatt68} 
Fattorini, H.O., Boundary control systems, {\em SIAM Journal on Control and Optimization}, {\bf 6}, 3, pp.\ 349--385, 1968.

\bibitem{Nikos} Kontaras, N., Heertjes,  M.F., and Zwart, H.J.  (2016). Continuous compliance compensation of position dependent flexible structures. {\it IFAC-PapersOnLine}, vol. 49, no. 13 , pp.\ 76-81.

\bibitem{Bauchau} O.A.\ Bauchau, J.I.\ Craig,  Chapter 16-Kirchhoff plate theory, Structural Analysis with applications to aerospace structures, pages 819-865. Springer, 2009.

\bibitem{Bogner} F.K.\ Bogner, R.L.\ Fox, and L.A.\ Schmit, The generation of interelementcompatible
stiffness and mass matrices by the use of interpolation formulae, in:
Procedings of the First Conference on Matrix Methods in Structural Mechanics,
vol. AFFDL-TR-66-80, Wright Patterson Air Force Base, Ohio, October
1966, pp. 397–443.

\bibitem{Brugnoli2} A.\ Brugnoli, D.\ Alazard, V.\ Pommier-Budinger, D.\ Matignon,  Port-Hamiltonian formulation and Symplectic discretization of Plate models Part II : Kirchhoff model for thin plate, {\em Applied Mathematical Modelling}, Vol.\ 75, 2019, pp.\ 961-981.

%\bibitem{Banks}.

%\bibitem{chill2019} R. Chill, L. Paunonen, D. Seifert, R. Stahn, Y. Tomilov (2019). Non-uniform stability of damped contraction semigroups. arXiv preprint. arXiv:1911.04804.

\bibitem{curtain2009} R.\ Curtain, K.\ Morris. Transfer functions of distributed parameter systems: A tutorial. {\it Automatica} Vol.\ 45, No.\ 5, 2009, pp.\ 1101-1116.

%\bibitem{haraux1989} A. Haraux. Une remarque sur la stabilisation de certains syst\`emes du deuxi\`eme ordre en temps. Portugaliae Mathematica 46 (1989):245-258.

\bibitem{rixen2014}  M. G{\'e}radin, D.J. Rixen. Mechanical vibrations: theory and application to structural dynamics. John Wiley \& Sons, 2014. 

\bibitem{Reddy} 
J.N. Reddy, {\em Theory and Analysis of Elastic Plates and Shells}. CRC press, 2006.

\bibitem{Lagnese}
J.E. Lagnese, {\em The Hilbert Uniqueness method: A retrospective}, Springer, Berlin, Heidelberg, 1991.

\bibitem{Renardy}
M. Renardy, R. Rogers, {\em An Introduction to Partial Differential Equations}. Texts in Applied Mathematics 13 (Second ed.). Springer-Verlag, New York, 2004.

\bibitem{Zienkiewicz}
O.C. Zienkiewicz, R.L. Taylor, J.Z. Zhu, {\em The Finite Element Method: its Basis and Fundamentals (Seventh Edition)}. Butterworth-Heinemann, Oxford, 2013. 

\end{thebibliography}
\end{document}